\documentclass[12pt,fleqn]{article}
\usepackage{amscd,amsmath,amssymb,amsthm}
\usepackage{geometry} 
\usepackage{thinsp}
\usepackage{tikz}
\usepackage[pdfpagemode=UseNone,pdfstartview=FitH]{hyperref}
\usepackage{doi}

\newcommand{\A}{{\mathcal A}}
\newcommand{\As}[1][]{A_\sigma #1}
\newcommand{\Bs}[1][]{B_\sigma #1}

\newcommand{\bgset}[1]{\big\{#1\big\}}

\newcommand{\dint}{\ds{\int}}
\newcommand{\ds}[1]{\displaystyle #1}
\newcommand{\eps}{\varepsilon}
\newcommand{\F}{{\mathcal F}}
\renewcommand{\H}{{\mathcal H}}
\newcommand{\M}{{\mathcal M}}
\newcommand{\N}{\mathbb N}
\newcommand{\norm}[2][]{\left\|#2\right\|_{#1}}
\renewcommand{\O}{\text{O}}
\renewcommand{\o}{\text{o}}
\newcommand{\PS}[1]{$(\text{PS})_{#1}$}
\newcommand{\pnorm}[2][]{\if #1'' \left|#2\right|_p \else \left|#2\right|_{#1} \fi}
\newcommand{\R}{\mathbb R}
\newcommand{\rad}{\text{rad}}
\newcommand{\restr}[2]{\left.#1\right|_{#2}}
\newcommand{\seq}[1]{\left(#1\right)}
\newcommand{\set}[1]{\left\{#1\right\}}
\newcommand{\wto}{\rightharpoonup}
\newcommand{\Z}{\mathbb Z}
\newcommand{\T}{{\mathcal T}}
\newcommand{\Tc}{{\mathcal T}_{\rm c}}
\newcommand{\Te}{{\mathcal T}_{\rm e}}
\newcommand{\Tsub}{{\mathcal T}_{\rm sub}}
\newcommand{\cE}{{\mathcal E}}

\newenvironment{enumroman}{\begin{enumerate}

}{\end{enumerate}}

\newtheorem{corollary}{Corollary}[section]
\newtheorem{lemma}[corollary]{Lemma}
\newtheorem{proposition}[corollary]{Proposition}
\newtheorem{theorem}[corollary]{Theorem}

\theoremstyle{definition}
\newtheorem{definition}[corollary]{Definition}

\theoremstyle{remark}

\newtheorem{remark}[corollary]{Remark}

\numberwithin{equation}{section}

\title{\bf A scaling-based approach to mixed fractional $p$-Laplacian equations in $\R^N$}

\author{%
\textbf{Kanishka Perera}\\[0.45em]
\parbox{0.96\textwidth}{%
\small\centering
Department of Mathematics\\
Florida Institute of Technology\\
150 W University Blvd, Melbourne, FL 32901-6975, USA\\[0.3em]
\textit{kperera@fit.edu}%
\bigskip \bigskip \bigskip
}\\[\bigskipamount]
\textbf{Caterina Sportelli}\\[0.45em]
\parbox[t]{0.43\textwidth}{%
\small\centering
Departamento de An\'alisis Matem\'atico\\
Universidad de Granada\\
Avenida de Fuente Nueva, s/n\\
18071 Granada, Spain\\[0.3em]
\textit{caterina.sp@ugr.es}%
}%
\hspace{0.03\textwidth}%
\parbox[t]{0.50\textwidth}{%
\small\centering
\mbox{Department of Mathematics and Statistics}\\
\mbox{The University of Western Australia}\\
35 Stirling Highway\\
Crawley, Perth, WA 6009, Australia\\[0.3em]
\textit{caterina.sportelli@uwa.edu.au}%
}%
}

\date{}

\begin{document}

\maketitle

\begin{abstract}
We study equations driven by a sum of two fractional $p$-Laplacians with different orders and exponents. A central contribution of this paper is a Poho\v{z}aev identity for finite-energy weak solutions, established without additional regularity assumptions by means of a new approach based on scaling and convexity defects. This identity, which is of independent interest, yields nonexistence results for critical problems. Furthermore, we establish weighted Sobolev embeddings, construct an unbounded sequence of positive eigenvalues without relying on the Poho\v{z}aev identity, and prove multiplicity results for subcritical and doubly critical problems using scaling-based variational methods.
\end{abstract}

\begin{center}
\begin{minipage}{12cm}
\footnotesize \tableofcontents
\end{minipage}
\end{center}

\newpage

\section{Introduction}

Elliptic problems involving nonlinear superposition operators of mixed fractional orders have been widely studied in recent years (see, e.g., \cite{MR3910033,MR4142328,MR4736013,MR4942302,MR4821750,KPCS} and the references therein). We recall that for $s \in (0,1)$ and $p \in (1,\infty)$, $(- \Delta)_p^s$ is the fractional $p$-Laplacian operator defined on smooth functions by
\[
(- \Delta)_p^s\, u(x) = 2\, \lim_{\eps \searrow 0} \int_{\R^N \setminus B_\eps(x)} \frac{|u(x) - u(y)|^{p-2}\, (u(x) - u(y))}{|x - y|^{N+sp}}\, dy, \quad x \in \R^N.
\]
Consider the model equation
\begin{equation} \label{1}
\A(u) := (- \Delta)_{p_1}^{s_1}\, u + (- \Delta)_{p_2}^{s_2}\, u = f(|x|,u) \quad \text{in } \R^N,
\end{equation}
where
\begin{gather}
\notag 0 < s_1 < s_2 < 1, \qquad 1 < p_2 < p_1 <\infty,\\[10pt]
\label{12} s_1 p_1 < s_2 p_2 < N < \frac{(s_2 - s_1)\, p_1 p_2}{p_1 - p_2} =: \alpha,
\end{gather}
and $f$ is a Carath\'{e}odory function on $[0,\infty) \times \R$ satisfying a suitable growth condition. When $p_1 \ne p_2$, the operator $\A$ is nonhomogeneous and therefore has no associated homogeneous eigenvalue problem. This creates difficulties in understanding the geometry of the energy functional associated with equation \eqref{1} and in applying standard variational methods. However, this operator has an important scaling property when $p_1 \ne p_2$. Let
\[
\beta = \frac{s_2 p_2 - s_1 p_1}{p_1 - p_2}
\]
and set
\[
u_t(x) = t^\beta\, u(tx), \quad x \in \R^N,\, t \ge 0.
\]
The two fractional $p$-Laplacians in $\A$ scale the same way under this scaling, resulting in
\begin{equation} \label{2}
\A(u_t)(x) = t^{\alpha - \beta}\, \A(u)(tx).
\end{equation}
This scaling property gives rise to a family of scaled eigenvalue problems with singular weights for the operator $\A$. In this paper we combine insights into the geometry of the energy functional gleaned from these scaled eigenvalue problems with variational methods for scaled functionals recently developed in Mercuri and Perera \cite{MR5043800} and Perera \cite{Pe25} to prove new existence and multiplicity results for equation \eqref{1}. A direct computation from the definitions of $\alpha$ and $\beta$ gives
\begin{equation} \label{alpha-beta}
\frac{\alpha}{p_i} = \beta + s_i, \quad i = 1, 2,
\end{equation}
an identity that will be used repeatedly.

Recall that the homogeneous fractional Sobolev space $\dot{W}^{s,p}(\R^N)$ is the completion of $C^\infty_c(\R^N)$ with respect to the Gagliardo-Slobodecki\u{\i} seminorm
\[
[u]_{\dot{W}^{s,p}(\R^N)} = \left(\int_{\R^N} \int_{\R^N} \frac{|u(x) - u(y)|^p}{|x - y|^{N+sp}}\, dx dy\right)^{1/p}.
\]
The natural space to look for solutions to equation \eqref{1} is the intersection space
\[
W = \dot{W}^{s_1,p_1}(\R^N) \cap \dot{W}^{s_2,p_2}(\R^N)
\]
endowed with the norm
\[
\norm{u} = [u]_1 + [u]_2, \qquad [u]_i := [u]_{\dot{W}^{s_i,p_i}(\R^N)}.
\]
We work in the subspace $W_\rad$ of radial functions in $W$.

\subsection{The embedding triangle}

For $i = 1,2$, $\dot{W}^{s_i,p_i}(\R^N)$ is continuously embedded in $L^{p_i^\ast}(\R^N)$ by the fractional Sobolev inequality
\[
\int_{\R^N} |u|^{p_i^\ast}\, dx \le C\, [u]_i^{p_i^\ast},
\]
where $p_i^\ast = Np_i/(N - s_i\thinthinspace p_i)$ is the fractional critical Sobolev exponent. It is also continuously embedded in the weighted Lebesgue space $L^{p_i}(\R^N;|x|^{-s_i p_i})$ by the fractional Hardy inequality
\[
\int_{\R^N} \frac{|u|^{p_i}}{|x|^{s_i p_i}}\, dx \le C\, [u]_i^{p_i}.
\]
A direct computation shows that
\begin{equation} \label{p-ast-order}
p_1^\ast < p_2^\ast \iff N\, (p_1 - p_2) < (s_2 - s_1)\, p_1 p_2 \iff N < \alpha,
\end{equation}
so $p_1^\ast < p_2^\ast$ precisely under \eqref{12}. For $q \ge 1$ and $b \ge 0$, let $L^q(\R^N;|x|^{-b})$ be the weighted Lebesgue space with the norm
\[
\norm[L^q(\R^N;|x|^{-b})]{u} := \left(\int_{\R^N} |x|^{-b}\, |u|^q\, dx\right)^{1/q}.
\]
We will show that the natural range of the embeddings of $W$ into weighted Lebesgue spaces is described by the triangle $\T$ in the $(q,b)$-plane with the vertices
\[
P(0,N), \qquad Q_1(p_1^\ast,0), \qquad Q_4(p_2^\ast,0)
\]
(see Figure \ref{fig:H}). The edges $P Q_1$ and $P Q_4$ lie on the lines
\[
b = N \left(1 - \frac{q}{p_1^\ast}\right) \quad \text{and} \quad b = N \left(1 - \frac{q}{p_2^\ast}\right),
\]
respectively, and they contain the points
\[
Q_2(p_1,s_1 p_1) \in P Q_1, \qquad Q_3(p_2,s_2 p_2) \in P Q_4,
\]
which correspond to the fractional Hardy inequalities in $\dot{W}^{s_1,p_1}(\R^N)$ and $\dot{W}^{s_2,p_2}(\R^N)$, respectively. We write $[A,B]$ for the closed segment joining the points $A$ and $B$. Set
\begin{align*}
\Tc & = \set{(q,b) : q \ge 1,\ b \ge 0,\ N\, (1 - q/p_1^\ast) < b < N\, (1 - q/p_2^\ast)},\\[5pt]
\Te & = \Tc \cup [Q_2,Q_1] \cup [Q_3,Q_4].
\end{align*}
Thus, $\Tc$ consists of the points of $\T$ with $q \ge 1$ that do not lie on the edges $P Q_1$ and $P Q_4$, and $\Te$ is obtained from $\set{(q,b) \in \T : q \ge 1}$ by removing the parts of the edges $P Q_1$ and $P Q_4$ that lie strictly between $P$ and the Hardy points $Q_2$ and $Q_3$. We will prove the following embedding result (see Theorem \ref{Theorem embedding}):
\begin{enumroman}
\item $W$ is continuously embedded in $L^q(\R^N;|x|^{-b})$ for all $(q,b) \in \Te$;
\item $W$ is compactly embedded in $L^q(\R^N;|x|^{-b})$ for all $(q,b) \in \Tc$ with $b > 0$, and $W_\rad$ is compactly embedded in $L^q(\R^N)$ for all $q \in (p_1^\ast\!,p_2^\ast)$ when $N \ge 2$;
\item the embedding of $W_\rad$ is not compact on $[Q_2,Q_1] \cup [Q_3,Q_4]$.
\end{enumroman}
In particular, $W$ is embedded in $L^q(\R^N)$ continuously for all $q \in [p_1^\ast\!,p_2^\ast]$, and $W_\rad$ is embedded in $L^q(\R^N)$ compactly for $q \in (p_1^\ast\!,p_2^\ast)$ when $N \ge 2$. The proof is elementary and does not use any pointwise (Strauss-type) estimates for radial functions: splitting the integral at $|x| = 1$, the region near the origin is controlled by the $\dot{W}^{s_2,p_2}$-part of the norm and the region near infinity by the $\dot{W}^{s_1,p_1}$-part. The convex hull of $Q_1$, $Q_2$, $Q_3$, and $Q_4$, given by
\begin{equation} \label{13}
\H = \set{\sum_{i=1}^4 \theta_i Q_i : \text{each } \theta_i \ge 0 \text{ and } \sum_{i=1}^4 \theta_i = 1},
\end{equation}
is the part of $\T$ that lies on or above the line through $Q_2$ and $Q_3$, which will turn out to be the scaling line (see Section \ref{Subsection regimes}), and $\T \setminus \H$ is the triangle $\Tsub$ with the vertices $P$, $Q_2$, and $Q_3$.

\begin{figure}[ht]
  \centering
  \begin{tikzpicture}[scale=1.0]
    \draw[->] (0,0) -- (7.2,0) node[right] {$q$};
    \draw[->] (0,0) -- (0,5.8) node[above] {$b$};
    \coordinate (P) at (0,5.2);
    \coordinate (Q3) at (1.5,3.98);
    \coordinate (Q2) at (2.6,1.82);
    \coordinate (Q1) at (4.0,0.0);
    \coordinate (Q4) at (6.4,0.0);
    \fill[gray!10] (P) -- (Q3) -- (Q2) -- cycle;
    \fill[gray!30] (Q3) -- (Q2) -- (Q1) -- (Q4) -- cycle;
    \draw[thick, dotted] (P) -- (Q3);
    \draw[thick, dotted] (P) -- (Q2);
    \draw[very thick] (Q2) -- (Q1);
    \draw[very thick] (Q3) -- (Q4);
    \draw (Q1) -- (Q4);
    \draw[dashed] (Q2) -- (Q3);
    \draw[densely dotted, black!60] (0.6,0) -- (0.6,4.71);
    \node[below] at (0.6,0) {\footnotesize $q = 1$};
    \fill (P) circle (1.7pt) node[left] {$P(0,N)$};
    \fill (Q1) circle (1.7pt) node[below] {$Q_1(p_1^\ast,0)$};
    \fill (Q2) circle (1.7pt) node[left=2pt] {$Q_2(p_1,s_1 p_1)$};
    \fill (Q3) circle (1.7pt) node[above right] {$Q_3(p_2,s_2 p_2)$};
    \fill (Q4) circle (1.7pt) node[below] {$Q_4(p_2^\ast,0)$};
    \node (sub) at (2.3,5.0) {\scriptsize subscaled};
    \draw[->, thin] (sub.south west) -- (1.15,3.95);
    \node at (4.75,0.55) {\footnotesize superscaled};
  \end{tikzpicture}
  \caption{The triangle $\T$ in the $(q,b)$-plane (schematic). $W$ is embedded in $L^q(\R^N;|x|^{-b})$ for $q \ge 1$ and $(q,b) \in \T$ off the dotted segments $[P,Q_2)$ and $[P,Q_3)$, compactly off the edges (for $W_\rad$ on the bottom edge). The dashed segment $Q_2 Q_3$ is the scaling line, which separates the subscaled region $\Tsub$ (light) from the superscaled region (dark). The thick segments $[Q_2,Q_1]$ and $[Q_3,Q_4]$ carry the critical exponents. The quadrilateral $\H$ is the dark region together with its boundary.}
  \label{fig:H}
\end{figure}

In view of these embeddings, a natural growth condition for the nonlinearity $f$ in equation \eqref{1} is
\begin{equation} \label{6}
|f(|x|,t)| \le \sum_{i=1}^n a_i\, |x|^{-b_i}\, |t|^{q_i - 1} \quad \forall (x,t) \in \R^N \times \R,
\end{equation}
where $a_i > 0$ are constants and each $(q_i,b_i) \in \Te$ with $q_i > 1$. Under this growth condition, a weak radial solution of equation \eqref{1} is a function $u \in W_\rad$ satisfying
\begin{multline*}
\int_{\R^N} \int_{\R^N} \left(\frac{|u(x) - u(y)|^{p_1 - 2}}{|x - y|^{N + s_1 p_1}} + \frac{|u(x) - u(y)|^{p_2 - 2}}{|x - y|^{N + s_2 p_2}}\right)\! (u(x) - u(y))(v(x) - v(y))\, dx dy\\[5pt]
= \int_{\R^N} f(|x|,u)\, v\, dx \quad \forall v \in W_\rad.
\end{multline*}
Weak radial solutions are the critical points of the $C^1$-functional
\[
E(u) = \int_{\R^N} \int_{\R^N} \left(\frac{1}{p_1}\, \frac{|u(x) - u(y)|^{p_1}}{|x - y|^{N + s_1 p_1}} + \frac{1}{p_2}\, \frac{|u(x) - u(y)|^{p_2}}{|x - y|^{N + s_2 p_2}}\right) dx dy - \int_{\R^N} F(|x|,u)\, dx
\]
defined on $W_\rad$, where $F(|x|,t) = \int_0^t f(|x|,\tau)\, d\tau$ is the primitive of $f$. By the principle of symmetric criticality of Palais \cite{MR547524}, weak radial solutions are also weak solutions, i.e., the above identity holds for all $v \in W$ (see Lemma \ref{Lemma symmetric criticality}).

\subsection{Nonlinear regimes} \label{Subsection regimes}

As in \cite{MR5043800}, the scaling property \eqref{2} of the operator $\A$ naturally leads to the following classification of different nonlinear regimes for $f$:
\begin{enumroman}
\item $f$ is subscaled if
    \[
    \lim_{|t| \to \infty}\, \frac{f(|t|^{-1},|t|^{\beta - 1}\, t)}{|t|^{\alpha - \beta - 1}\, t} = 0;
    \]
\item $f$ is asymptotically scaled if
    \[
    \lim_{|t| \to \infty}\, \frac{f(|t|^{-1},|t|^{\beta - 1}\, t)}{|t|^{\alpha - \beta - 1}\, t} = \lambda \in (0,\infty);
    \]
\item $f$ is superscaled if
    \[
    \lim_{|t| \to \infty}\, \frac{f(|t|^{-1},|t|^{\beta - 1}\, t)}{|t|^{\alpha - \beta - 1}\, t} = \infty.
    \]
\end{enumroman}
We will see that the energy functional $E$ has distinct geometries in these three regimes, which leads to different types of existence and multiplicity results for equation \eqref{1}.

First we consider the pure weighted power equation
\begin{equation} \label{7}
(- \Delta)_{p_1}^{s_1}\, u + (- \Delta)_{p_2}^{s_2}\, u = \lambda\, |x|^{-b}\, |u|^{q-2}\, u \quad \text{in } \R^N,
\end{equation}
where $(q,b) \in \Te$. For the nonlinearity $f(|x|,t) = \lambda\, |x|^{-b}\, |t|^{q-2}\, t$, we have
\[
\frac{f(|t|^{-1},|t|^{\beta - 1}\, t)}{|t|^{\alpha - \beta - 1}\, t} = \lambda\, |t|^{b - \alpha + \beta q}.
\]
By \eqref{alpha-beta}, the line $b = \alpha - \beta q$ passes through $Q_2$ and $Q_3$; we call it the scaling line. Since $\alpha > N$, the vertex $P(0,N)$ of $\T$ lies strictly below the scaling line, while $Q_1$ and $Q_4$ lie on or above it. Hence $b < \alpha - \beta q$ in $\Tsub \setminus [Q_2,Q_3]$, where $f$ is subscaled, and $b \ge \alpha - \beta q$ in $\H$, with equality exactly on the edge $Q_2 Q_3$. If $b > \alpha - \beta q$, then $f$ is superscaled and we will see that equation \eqref{7} has infinitely many solutions with positive energy for all $\lambda > 0$ in the subcritical case $(q,b) \notin [Q_2,Q_1] \cup [Q_3,Q_4]$, under conditions made precise in Theorem \ref{Theorem 2}. If $b = \alpha - \beta q$, then $f$ is asymptotically scaled and equation \eqref{7} is a scaled eigenvalue problem in the sense of \cite{MR5043800}, i.e., if $u$ is a solution of this equation, then the entire $1$-parameter family of functions $\set{u_t : t > 0}$ are solutions. When $(q,b) \in Q_2 Q_3 \setminus \set{Q_2,Q_3}$, the minimax scheme of \cite{MR5043800} based on the $\Z_2$-cohomological index of Fadell and Rabinowitz \cite{MR0478189} gives us an unbounded sequence of positive eigenvalues with radial eigenfunctions (see Theorem \ref{Theorem 1}).

\subsection{Poho\v{z}aev identities} \label{Subsection Pohozaev review}

The abstract theory of \cite{MR5043800,Pe25} requires that {\em every} solution of the relevant equation satisfies a Poho\v{z}aev-type identity (see assumptions $(A_{11})$ and $(f_4)$ in Section \ref{Section abstract}). For problem \eqref{7} this is the identity
\begin{equation} \label{Pohozaev intro}
\sigma \sum_{i=1}^2 \frac{1}{p_i} \int_{\R^N} \int_{\R^N} \frac{|u(x) - u(y)|^{p_i}}{|x - y|^{N + s_i p_i}}\, dx dy = \frac{\rho \lambda}{q} \int_{\R^N} |x|^{-b}\, |u|^q\, dx, \qquad \sigma = \alpha - N, \quad \rho = b - N + \beta q,
\end{equation}
which expresses the vanishing at $t = 1$ of the derivative of $E(u_t)$. Establishing \eqref{Pohozaev intro} for all finite-energy solutions turns out to be the main analytical difficulty of the present work, and we now review what is known.

For the Laplacian, such identities go back to Poho\v{z}aev \cite{MR192184}, and for general quasilinear operators to Pucci and Serrin \cite{MR855181}, who assumed $C^2$-regularity. For the $p$-Laplacian, Guedda and V\'{e}ron \cite{MR1009077} treated solutions in $W^{1,p}_0(\Omega) \cap L^\infty(\Omega)$, Degiovanni, Musesti, and Squassina \cite{MR2018671} showed that $C^1$-regularity suffices, and Liu and Liu \cite{MR2784366} obtained the identity in $\R^N$ for solutions in $W^{1,p}(\R^N)$ by testing the equation with difference quotients. For the fractional Laplacian, the identity was proved by Ros-Oton and Serra \cite{MR3211861} in bounded domains and by Chang and Wang \cite{MR3007900} in $\R^N$; see also Djitte, Fall, and Weth \cite{MR4685066} for a generalized fractional Poho\v{z}aev identity. For the fractional $p$-Laplacian with $p \ne 2$ the problem is considerably more delicate: there is no extension method, and the available regularity theory does not reach $C^1$ in general. In $\R^N$, Ambrosio \cite{MR4761133} proved the identity for solutions of class $C^{1,1}(\R^N)$, and Anthal and Garain \cite{AnGa25} adapted the difference-quotient method of \cite{MR2784366} to solutions in $W^{1,p}(\R^N) \cap W^{s,p}(\R^N)$.

For operators that combine two different orders of differentiation, the literature is much more limited. In the mixed local-nonlocal case, Biswas \cite{MR4992640} proved a Poho\v{z}aev identity for the semilinear operator $- \Delta + a\, (- \Delta)^s$, $a \ge 0$, in bounded $C^2$-domains (see also \cite{MR4597987,MR4785698} for related problems), and Anthal and Garain \cite{AnGa25} proved one in $\R^N$ for the quasilinear operator $- H_p + \gamma\, (- \Delta)_p^s$, where $H_p$ is an anisotropic $p$-Laplacian and {\em both operators have the same exponent $p$}, for solutions in $W^{1,p}(\R^N) \cap W^{s,p}(\R^N) \cap C^{0,l}(\R^N)$ with $l > s$ and nonlinearities of $(p - 1)$-growth (or bounded solutions); this identity was used by Biswas, Chakraborty, and Das \cite{BiChDa25} to obtain nonexistence results for critical problems for $- \Delta_p + (- \Delta)_p^s$ in $\R^N$. In the linear case $p_1 = p_2 = 2$, the operator $(- \Delta)^{s_1} + (- \Delta)^{s_2}$ is a Fourier multiplier and a Poho\v{z}aev identity in $\R^N$ is available (see, e.g., Chergui, Gou, and Hajaiej \cite[Lemma 2.2]{MR4626463}), where it is used in the study of normalized solutions; see also Luo and Hajaiej \cite{MR4445936}. To the best of our knowledge, however, no Poho\v{z}aev identity is available in the literature for a sum of two fractional $p$-Laplacians when $(p_1,p_2) \ne (2,2)$, and in particular none for the nonhomogeneous operator $\A$ with $p_1 \ne p_2$.

In this paper we prove such an identity (see Theorem \ref{Theorem Pohozaev}): if $(q,b) \in \H$ and
\begin{equation} \label{C}
2b > \max \set{q,2}, \tag{C}
\end{equation}
then {\em every} finite-energy radial weak solution of equation \eqref{7} satisfies \eqref{Pohozaev intro}. No regularity or boundedness of the solution is assumed, and $(q,b)$ is allowed to lie on the critical edges $[Q_2,Q_1]$ and $[Q_3,Q_4]$. The proof is not based on the Euler multiplier $x \cdot \nabla u$, which is not available for energy solutions, but on testing the equation with $u_t - u$, where $u_t$ is the scaling of the problem. Two observations make this work. First, both Gagliardo energies and the weighted potential are {\em exactly} covariant under the scaling $u \mapsto u_t$, so the expansion of the energy along the curve $t \mapsto u_t$ is explicit up to second-order remainders. Second, these remainders are convexity defects. The defects of the operator energies are nonnegative and need not be estimated at all, while that of the weighted potential is superlinearly small in $t - 1$ by an estimate for the dilation modulus of fractional Sobolev functions in Hardy-weighted norms (Lemma \ref{Lemma modulus}). The convexity-defect mechanism is inspired by the difference-quotient proofs of \cite{MR2784366,AnGa25}, which test with $\sum_j x_j\, D_j u$ and squeeze the defects through the equation; replacing the coordinate translations by the dilation group of the problem eliminates all commutators, cutoffs, and principal-value manipulations, and with them all regularity requirements. The weight $|x|^{-b}$, usually a nuisance, is here the resource: condition \eqref{C} fails at $b = 0$, consistently with the fact that for unweighted problems the identities of \cite{MR2784366,AnGa25} require $\nabla u \in L^p$ and that of \cite{MR4761133} requires $C^{1,1}$-regularity.

The Poho\v{z}aev identity is used in this paper only where the scaled-manifold methods of \cite{MR5043800,Pe25} are indispensable. In the superscaled case with $q > p_1$, as well as for the Brezis-Nirenberg type critical problems considered in Section \ref{Subsection BN}, we give direct proofs of the boundedness of Palais-Smale sequences, and no Poho\v{z}aev identity is needed.

\section{Statements of main theorems}

\subsection{A Poho\v{z}aev identity}

Let $\sigma = \alpha - N$ and $\rho = b - N + \beta q$ be as in \eqref{Pohozaev intro}. Our first main result is the following Poho\v{z}aev-type identity for the pure weighted power equation \eqref{7}.

\begin{theorem} \label{Theorem Pohozaev}
Let $(q,b) \in \H$ satisfy \eqref{C}, and let $\lambda \in \R$. Then every weak radial solution $u \in W_\rad$ of equation \eqref{7} satisfies the Poho\v{z}aev identity \eqref{Pohozaev intro}.
\end{theorem}

\begin{remark} \label{Remark Pohozaev region}
On the scaling edge $Q_2 Q_3$ we have $b = \alpha - \beta q$, which decreases from $s_2 p_2$ at $q = p_2$ to $s_1 p_1$ at $q = p_1$, while $\max \set{q,2}$ is nondecreasing in $q$. Hence the set of $q \in (p_2,p_1)$ for which \eqref{C} holds on this edge is an interval $(p_2,\bar{q})$, which is nonempty if and only if $2 s_2 p_2 > \max \set{p_2,2}$ and equal to $(p_2,p_1)$ if $2 s_1 p_1 \ge \max \set{p_1,2}$. For example, $N = 3$, $s_1 = 0.6$, $s_2 = 0.9$, $p_1 = 2.2$, and $p_2 = 2$ satisfy \eqref{12} with $\alpha = 6.6$, and $2 s_1 p_1 = 2.64 > p_1$, so \eqref{C} holds on the whole edge $Q_2 Q_3$ in this case.
\end{remark}

\begin{remark}
The proof of Theorem \ref{Theorem Pohozaev} also applies to non-radial solutions $u \in W$, i.e., functions satisfying the weak formulation of equation \eqref{7} for all $v \in W$, since $u_t - u \in W$ for all $t > 0$. It uses only three properties of the energy: it is a sum of convex functionals that are exactly covariant under a common one-parameter scaling, the potential is covariant under the same scaling, and the dilation modulus estimate of Corollary \ref{Corollary scaled modulus}. Consequently, the same argument applies to superpositions of more than two fractional $p$-Laplacians admitting a common scaling, and to the mixed local-nonlocal case $s_2 = 1$, where the corresponding modulus estimate for $\dot{W}^{1,p_2}(\R^N)$ follows from the Hardy inequality.
\end{remark}

\subsection{Pure weighted power case: \texorpdfstring{$f(|x|,t) = \lambda\, |x|^{-b}\, |t|^{q-2}\, t$}{}}

Consider the pure weighted power equation
\begin{equation} \label{15}
(- \Delta)_{p_1}^{s_1}\, u + (- \Delta)_{p_2}^{s_2}\, u = \lambda\, |x|^{-b}\, |u|^{q-2}\, u \quad \text{in } \R^N,
\end{equation}
where $(q,b) \in \H$. The energy functional associated with this equation is
\[
E(u) = \int_{\R^N} \int_{\R^N} \left(\frac{1}{p_1}\, \frac{|u(x) - u(y)|^{p_1}}{|x - y|^{N + s_1 p_1}} + \frac{1}{p_2}\, \frac{|u(x) - u(y)|^{p_2}}{|x - y|^{N + s_2 p_2}}\right) dx dy - \lambda \int_{\R^N} \frac{1}{q}\, |x|^{-b}\, |u|^q\, dx
\]
defined on $W_\rad$.

As we have noted in the introduction, here
\[
\frac{f(|t|^{-1},|t|^{\beta - 1}\, t)}{|t|^{\alpha - \beta - 1}\, t} = \lambda\, |t|^{b - \alpha + \beta q}
\]
and $b - \alpha + \beta q \ge 0$ for $(q,b) \in \H$ since the edge $Q_2 Q_3$ of $\H$ lies on the scaling line $b = \alpha - \beta q$. This together with Theorem \ref{Theorem embedding} leads to the following classification of different nonlinear regimes:
\begin{enumroman}
\item if $(q,b) \in Q_2 Q_3 \setminus \set{Q_2,Q_3}$, then equation \eqref{15} is subcritical and asymptotically scaled;
\item if $(q,b) \in \H \setminus ([Q_2,Q_1] \cup [Q_2,Q_3] \cup [Q_3,Q_4])$, then equation \eqref{15} is subcritical and superscaled;
\item if $(q,b) \in [Q_2,Q_1] \cup [Q_3,Q_4]$, then equation \eqref{15} is critical.
\end{enumroman}
We note that the equations of the edges $Q_1 Q_2$ and $Q_3 Q_4$ are $b = N\, (1 - q/p_1^\ast)$ and $b = N\, (1 - q/p_2^\ast)$, respectively.

\subsubsection{Subcritical and asymptotically scaled case: \texorpdfstring{$p_2 < q < p_1$ \normalfont{and} $b = \alpha - \beta q$}{}}

In this case equation \eqref{15} is a scaled eigenvalue problem (see Section \ref{Scaled eigenvalue problems}) and we have the following theorem.

\begin{theorem} \label{Theorem 1}
If $p_2 < q < p_1$ and $b = \alpha - \beta q$, then the scaled eigenvalue problem \eqref{15} has a sequence of positive eigenvalues $\widehat{\lambda}_k(q) \to \infty$ with radial eigenfunctions.
\end{theorem}

The eigenvalues $\widehat{\lambda}_k(q)$ are obtained from the minimax values $\lambda_k(q)$ of Section \ref{Scaled eigenvalue problems}, which are based on the $\Z_2$-cohomological index, and satisfy
\[
\frac{p_2}{q}\, \lambda_k(q) \le \widehat{\lambda}_k(q) \le \frac{p_1}{q}\, \lambda_k(q)
\]
(see Section \ref{Subsection proof Theorem 1}). If, in addition, $2b > \max \set{q,2}$, then $\widehat{\lambda}_k(q) = \lambda_k(q)$, so $\lambda_k(q) \nearrow \infty$ is a sequence of eigenvalues. By Remark \ref{Remark Pohozaev region}, this condition is satisfied on a nonempty subinterval $(p_2,\bar{q})$ of $(p_2,p_1)$ when $2 s_2 p_2 > \max \set{p_2,2}$, and on all of $(p_2,p_1)$ when $2 s_1 p_1 \ge \max \set{p_1,2}$.

In fact, the equality $\widehat{\lambda}_1(q)=\lambda_1(q)$ holds without assuming \eqref{C}. Moreover, $\lambda_1(q)$ has a nonnegative radial eigenfunction $u_1$ satisfying
\[
\frac{[u_1]_1^{p_1}}{[u_1]_2^{p_2}} = \frac{p_1(q-p_2)}{p_2(p_1-q)}.
\]
These assertions follow from the variational characterization of $\lambda_1(q)$ and are proved at the end of Section \ref{Subsection proof Theorem 1}.

\subsubsection{Subcritical and superscaled case: \texorpdfstring{$N\, (1 - q/p_1^\ast) < b < N\, (1 - q/p_2^\ast)$ \normalfont{and} $b > \alpha - \beta q$}{}}

In this case we have the following theorem.

\begin{theorem} \label{Theorem 2}
Let $(q,b) \in \H$ satisfy $N\, (1 - q/p_1^\ast) < b < N\, (1 - q/p_2^\ast)$ and $b > \alpha - \beta q$, let $\lambda > 0$, assume that $N \ge 2$ if $b = 0$, and assume that either
\begin{enumroman}
\item $q > p_1$, or
\item $2b > \max \set{q,2}$.
\end{enumroman}
Then equation \eqref{15} has an infinite sequence of nontrivial radial solutions $\seq{u_k}$ such that $E(u_k) > 0$ and $E(u_k) \nearrow \infty$.
\end{theorem}

In case $(i)$, which contains all superscaled points of $\H$ with $q > p_1$, the proof is direct: Palais-Smale sequences are bounded by an argument of Ambrosetti-Rabinowitz type, and no Poho\v{z}aev identity is needed. In case $(ii)$ we use the constrained approach of \cite{Pe25} together with Theorem \ref{Theorem Pohozaev}. Since $b > \alpha - \beta q \ge s_1 p_1$ when $p_2 < q \le p_1$, the condition $2 s_1 p_1 \ge \max \set{p_1,2}$ ensures that Theorem \ref{Theorem 2} covers the whole subcritical superscaled region.

In particular, we have the following corollary for the subcritical pure power equation
\begin{equation} \label{27}
(- \Delta)_{p_1}^{s_1}\, u + (- \Delta)_{p_2}^{s_2}\, u = \lambda\, |u|^{q-2}\, u \quad \text{in } \R^N
\end{equation}
with the associated energy functional
\[
E(u) = \int_{\R^N} \int_{\R^N} \left(\frac{1}{p_1}\, \frac{|u(x) - u(y)|^{p_1}}{|x - y|^{N + s_1 p_1}} + \frac{1}{p_2}\, \frac{|u(x) - u(y)|^{p_2}}{|x - y|^{N + s_2 p_2}}\right) dx dy - \lambda \int_{\R^N} \frac{1}{q}\, |u|^q\, dx
\]
defined on $W_\rad$.

\begin{corollary} \label{Corollary pure power}
If $N \ge 2$, $p_1^\ast < q < p_2^\ast$, and $\lambda > 0$, then equation \eqref{27} has an infinite sequence of nontrivial radial solutions $\seq{u_k}$ such that $E(u_k) > 0$ and $E(u_k) \nearrow \infty$.
\end{corollary}

\subsubsection{Critical case: \texorpdfstring{$b = N\, (1 - q/p_1^\ast)$ \normalfont{or} $b = N\, (1 - q/p_2^\ast)$}{}}

Theorem \ref{Theorem Pohozaev} holds on the critical edges of $\H$, and gives the following nonexistence result.

\begin{theorem} \label{Theorem nonexistence}
Let $(q,b) \in [Q_2,Q_1] \cup [Q_3,Q_4]$ satisfy $2b > \max \set{q,2}$, and let $\lambda \in \R$. Then equation \eqref{15} has no nontrivial weak radial solution $u \in W_\rad$.
\end{theorem}

\begin{remark}
Condition \eqref{C} fails at the Sobolev points $Q_1$ and $Q_4$, where $b = 0$. Whether the unweighted pure critical power equations
\[
(- \Delta)_{p_1}^{s_1}\, u + (- \Delta)_{p_2}^{s_2}\, u = |u|^{p_i^\ast - 2}\, u \quad \text{in } \R^N, \quad i = 1, 2
\]
have nontrivial radial solutions is an open problem. For the double phase operator, the corresponding nonexistence results were obtained in \cite{PeWi} using a Poho\v{z}aev identity valid for all finite-energy solutions.
\end{remark}

\subsection{Brezis-Nirenberg type critical problems} \label{Subsection BN}

In this subsection we consider the doubly critical problem
\begin{equation} \label{BN}
(- \Delta)_{p_1}^{s_1}\, u + (- \Delta)_{p_2}^{s_2}\, u = \mu\, |u|^{q-2}\, u + \mu_1\, |u|^{p_1^\ast - 2}\, u + \mu_2\, |u|^{p_2^\ast - 2}\, u \quad \text{in } \R^N,
\end{equation}
where
\begin{equation} \label{BN hypotheses}
N \ge 2, \qquad p_1^\ast < q < p_2^\ast, \qquad \mu > 0, \qquad \mu_1, \mu_2 \ge 0, \qquad \mu_1 + \mu_2 > 0.
\end{equation}
By \eqref{p-ast-order}, the interval $(p_1^\ast,p_2^\ast)$ is nonempty precisely under \eqref{12}, so problem \eqref{BN} is naturally posed on the whole parameter domain \eqref{12}. The energy functional associated with this equation is
\[
E_\mu(u) = \sum_{i=1}^2 \frac{1}{p_i}\, [u]_i^{p_i} - \frac{\mu}{q} \int_{\R^N} |u|^q\, dx - \sum_{i=1}^2 \frac{\mu_i}{p_i^\ast} \int_{\R^N} |u|^{p_i^\ast}\, dx, \quad u \in W_\rad.
\]
Let
\[
S_i = \inf_{u \in \dot{W}^{s_i,p_i}(\R^N) \setminus \set{0}}\, \frac{[u]_i^{p_i}}{\left(\dint_{\R^N} |u|^{p_i^\ast}\, dx\right)^{p_i/p_i^\ast}}, \quad i = 1, 2
\]
be the best constants in the fractional Sobolev inequalities, set
\begin{equation} \label{beta-i}
\beta_i = \frac{s_i}{N}\, S_i^{N/s_i p_i}\, \mu_i^{- (N - s_i p_i)/s_i p_i}, \quad i = 1, 2,
\end{equation}
with the convention that $\beta_i = \infty$ if $\mu_i = 0$, and let
\[
\beta^\ast = \min \set{\beta_1,\beta_2}.
\]
Finally, fix a nested sequence $V_1 \subset V_2 \subset \cdots$ of subspaces of radial functions in $C^\infty_c(\R^N)$ with $\dim V_m = m$, and set
\begin{equation} \label{delta-m}
\delta_m = \min_{u \in V_m,\, \norm{u} = 1}\, \int_{\R^N} |u|^q\, dx > 0.
\end{equation}

\begin{theorem} \label{Theorem BN}
Assume \eqref{12} and \eqref{BN hypotheses}. Then for each $m \in \N$, problem \eqref{BN} has $m$ distinct pairs of nontrivial radial solutions $\pm\, u_1^\mu,\dots,\pm\, u_m^\mu$ with $0 < E_\mu(u_k^\mu) < \beta^\ast$ for all
\begin{equation} \label{nu-m}
\mu > \nu_m := \frac{q}{\delta_m}\, \sup_{\tau > 0} \left[\frac{1}{p_1\, \tau^{q - p_1}} + \frac{1}{p_2\, \tau^{q - p_2}} - \frac{\beta^\ast}{\tau^q}\right].
\end{equation}
In particular, the number of solutions of problem \eqref{BN} goes to infinity as $\mu \to \infty$.
\end{theorem}

\begin{remark}
The supremum in \eqref{nu-m} is finite and positive since the expression in brackets tends to $- \infty$ as $\tau \to 0$, tends to $0$ as $\tau \to \infty$, and is positive for large $\tau$. Moreover, $\delta_m \to 0$ as $m \to \infty$: otherwise the $L^q$-norm and the norm of $W$ would be equivalent on the infinite-dimensional space $\bigcup_m V_m$, contradicting the compactness of the embedding $W_\rad \hookrightarrow L^q(\R^N)$. Hence $\nu_m \to \infty$, and $1/\delta_m$ plays here the role played by the cohomological-index eigenvalues in the corresponding results for $p$-Laplacian and $(p,q)$-Laplacian problems in bounded domains in Perera \cite{MR4999814}.
\end{remark}

\begin{remark}
The proof of Theorem \ref{Theorem BN} uses neither a Poho\v{z}aev identity nor a profile decomposition, and nothing is assumed or used about the extremals of the fractional Sobolev inequalities: the compactness threshold $\beta^\ast$ is obtained by the elementary argument of \cite{PeWi}, combined with the inequality $E_\mu(u) \ge 0$ for all critical points $u$, which here follows directly from the ordering $p_2 < p_1 < p_1^\ast < q < p_2^\ast$ of the exponents (see Proposition \ref{Proposition PS}). The threshold $\beta^\ast$ does not depend on $\mu$, while the minimax levels used in the proof go to zero as $\mu \to \infty$; this is the mechanism behind the theorem, as in \cite{MR4999814}.
\end{remark}

\section{Scaling-based methods in critical point theory} \label{Section abstract}

In this section we recall from \cite{MR5043800} the construction of eigenvalues based on the cohomological index for scaled eigenvalue problems, and some abstract critical point theorems from \cite{MR4999814,Pe25}.

\subsection{Scalings}

Let $W$ be a reflexive Banach space. The following notion of a scaling in $W$ was introduced in \cite{MR5043800}.

\begin{definition}[{\cite[Definition 1.1]{MR5043800}}]
A scaling in $W$ is a continuous mapping
\[
W \times [0,\infty) \to W, \quad (u,t) \mapsto u_t
\]
satisfying
\begin{enumerate}
\item[$(A_1)$] $(u_{t_1})_{t_2} = u_{t_1 t_2}$ for all $u \in W$ and $t_1, t_2 \ge 0$;
\item[$(A_2)$] $(\tau u)_t = \tau u_t$ for all $u \in W$, $\tau \in \R$, and $t \ge 0$;
\item[$(A_3)$] $u_0 = 0$ and $u_1 = u$ for all $u \in W$;
\item[$(A_4)$] $u_t$ is bounded on bounded sets in $W \times [0,\infty)$;
\item[$(A_5)$] there exists $\sigma > 0$ such that $\norm{u_t} = \O(t^\sigma)$ as $t \to \infty$, uniformly in $u$ on bounded sets.
\end{enumerate}
\end{definition}

\subsection{Scaled operators}

Denote by $W^\ast$ the dual of $W$. Recall that $q \in C(W,W^\ast)$ is a potential operator if there is a functional $Q \in C^1(W,\R)$, called a potential for $q$, with Fr\'{e}chet derivative $Q' = q$. By replacing $Q$ with $Q - Q(0)$ if necessary, we may assume that $Q(0) = 0$.

\begin{definition}[{\cite[Definition 2.1]{MR5043800}}]
A scaled operator is an odd potential operator $\As \in C(W,W^\ast)$ that maps bounded sets into bounded sets and satisfies
\[
\As(u_t) v_t = t^\sigma \As(u) v \quad \forall u, v \in W,\, t \ge 0.
\]
\end{definition}

\subsection{Scaled eigenvalue problems} \label{Scaled eigenvalue problems}

Let $\As$ and $\Bs$ be scaled operators satisfying
\begin{enumerate}
\item[$(A_6)$] $\As(u) u > 0$ for all $u \in W \setminus \set{0}$;
\item[$(A_7)$] every sequence $\seq{u_j}$ in $W$ such that $u_j \wto u$ and $\As(u_j)(u_j - u) \to 0$ has a subsequence that converges strongly to $u$;
\item[$(A_8)$] $\Bs(u) u > 0$ for all $u \in W \setminus \set{0}$;
\item[$(A_9)$] if $u_j \wto u$ in $W$, then $\Bs(u_j) \to \Bs(u)$ in $W^\ast$.
\end{enumerate}
The scaled eigenvalue problem
\begin{equation} \label{8}
\As(u) = \lambda \Bs(u) \quad \text{in } W^\ast
\end{equation}
was studied in \cite{MR5043800}. We say that $\lambda \in \R$ is an eigenvalue of this problem if there is a $u \in W \setminus \set{0}$, called an eigenfunction associated with $\lambda$, satisfying \eqref{8}. Then $u_t$ is also an eigenfunction associated with $\lambda$ for any $t > 0$ since
\[
\As(u_t) v = \As(u_t) (v_{t^{-1}})_t = t^\sigma \As(u) v_{t^{-1}} = \lambda t^\sigma \Bs(u) v_{t^{-1}} = \lambda \Bs(u_t) (v_{t^{-1}})_t = \lambda \Bs(u_t) v
\]
for all $v \in W$.

The potentials
\begin{equation} \label{9}
I_\sigma(u) = \int_0^1 \As(\tau u) u\, d\tau, \quad J_\sigma(u) = \int_0^1 \Bs(\tau u) u\, d\tau, \quad u \in W
\end{equation}
of $\As$ and $\Bs$, respectively, are even, \hspace{-5pt} bounded \hspace{-5pt} on \hspace{-5pt} bounded sets, and have the scaling property
\[
I_\sigma(u_t) = t^\sigma I_\sigma(u), \quad J_\sigma(u_t) = t^\sigma J_\sigma(u) \quad \forall u \in W,\, t \ge 0
\]
(see \cite[Proposition 2.2]{MR5043800}). By \eqref{9}, $(A_6)$, and $(A_8)$,
\begin{equation} \label{16}
I_\sigma(u) > 0, \quad J_\sigma(u) > 0 \quad \forall u \in W \setminus \set{0}.
\end{equation}
It follows from $(A_9)$ that if $u_j \wto u$ in $W$, then $J_\sigma(u_j) \to J_\sigma(u)$ (see {\cite[Proposition 2.4]{MR5043800}}). We assume that $I_\sigma$ and $J_\sigma$ satisfy
\begin{enumerate}
\item[$(A_{10})$] $I_\sigma$ is coercive, i.e., $I_\sigma(u) \to \infty$ as $\norm{u} \to \infty$;
\item[$(A_{11})$] every solution of problem \eqref{8} satisfies the Poho\v{z}aev-type identity $I_\sigma(u) = \lambda J_\sigma(u)$.
\end{enumerate}

We have $I_\sigma'(0) = \As(0) = 0$ since $\As$ is odd, so the origin is a critical point of $I_\sigma$. It is the only critical point of $I_\sigma$ since
\[
I_\sigma'(u) u = \As(u) u > 0 \quad \forall u \in W \setminus \set{0}
\]
by $(A_6)$. So $I_\sigma(0) = 0$ is the only critical value of $I_\sigma$, and hence it follows from the implicit function theorem that
\[
\M = \bgset{u \in W : I_\sigma(u) = 1}
\]
is a $C^1$-Finsler manifold. Since $I_\sigma$ is continuous, even, and coercive, $\M$ is complete, symmetric, and bounded. Let
\[
\Psi(u) = \frac{1}{J_\sigma(u)}, \quad u \in W \setminus \set{0}
\]
and let $\widetilde{\Psi} = \restr{\Psi}{\M}$. Then eigenvalues of problem \eqref{8} coincide with critical values of $\widetilde{\Psi}$ (see \cite[Proposition 2.5]{MR5043800}).

Denote by $\F$ the class of symmetric subsets of $\M$ and by $i(M)$ the cohomological index of $M \in \F$. For $k \ge 1$, let $\F_k = \bgset{M \in \F : i(M) \ge k}$ and set
\[
\lambda_k := \inf_{M \in \F_k}\, \sup_{u \in M}\, \widetilde{\Psi}(u).
\]
The following theorem was proved in \cite{MR5043800}.

\begin{theorem}[{\cite[Theorem 2.10]{MR5043800}}] \label{Theorem 3}
Assume $(A_1)$--$(A_{11})$. Then $\lambda_k \nearrow \infty$ is a sequence of eigenvalues of problem \eqref{8}.
\begin{enumroman}
\item The first eigenvalue is given by
    \[
    \lambda_1 = \min_{u \in \M}\, \widetilde{\Psi}(u) > 0.
    \]
\item If $\lambda_k = \dotsb = \lambda_{k+m-1} = \lambda$ and $E_\lambda$ is the set of eigenfunctions associated with $\lambda$ that lie on $\M$, then $i(E_\lambda) \ge m$.
\item If $\lambda_k < \lambda < \lambda_{k+1}$, then
    \[
    i(\widetilde{\Psi}^{\lambda_k}) = i(\M \setminus \widetilde{\Psi}_\lambda) = i(\widetilde{\Psi}^\lambda) = i(\M \setminus \widetilde{\Psi}_{\lambda_{k+1}}) = k,
    \]
    where $\widetilde{\Psi}^a = \bgset{u \in \M : \widetilde{\Psi}(u) \le a}$ and $\widetilde{\Psi}_a = \bgset{u \in \M : \widetilde{\Psi}(u) \ge a}$ for $a \in \R$.
\end{enumroman}
\end{theorem}

\subsection{Multiplicity results for even functionals}

\subsubsection{Unconstrained case}

Let $S = \set{u \in W : \norm{u} = 1}$ be the unit sphere in $W$. The following abstract multiplicity theorem was proved in Perera \cite{MR4999814}. It will be used in the proof of Theorem \ref{Theorem BN}.

\begin{theorem}[{\cite[Theorem 2.1]{MR4999814}}] \label{Theorem pseudo-index}
Let $E \in C^1(W,\R)$ be an even functional and assume that there exists $c^\ast > 0$ such that $E$ satisfies the {\em \PS{c}} condition for all $c \in (0,c^\ast)$. Let $C$ be a compact symmetric subset of $S$ with $i(C) = m \ge 1$. Assume that the origin is a strict local minimizer of $E$ and that there exists $R > 0$ such that
\[
\sup_{u \in A}\, E(u) \le 0, \qquad \sup_{u \in X}\, E(u) < c^\ast,
\]
where $A = \set{Ru : u \in C}$ and $X = \set{tu : u \in A,\, 0 \le t \le 1}$. Then $E$ has $m$ distinct pairs of critical points at levels in $(0,c^\ast)$.
\end{theorem}

\subsubsection{Constrained case}

Let $\As$ be as in the last section. We consider the nonlinear operator equation
\begin{equation} \label{18}
\As(u) = f(u) \quad \text{in } W^\ast,
\end{equation}
where $f \in C(W,W^\ast)$ is an odd potential operator that maps bounded sets into bounded sets and satisfies
\begin{enumerate}
\item[$(f_1)$] there exists $\rho \ne \sigma$ such that $f(u_t) v_t = t^\rho f(u) v$ for all $u, v \in W$ and $t \ge 0$;
\item[$(f_2)$] $f(u) u > 0$ for all $u \in W \setminus \set{0}$;
\item[$(f_3)$] if $u_j \wto u$ in $W$, then $f(u_j) \to f(u)$ in $W^\ast$;
\item[$(f_4)$] for any $\mu \in \R$, every solution of the equation $\As(u) = \mu f(u)$ satisfies the Poho\v{z}aev-type identity $\sigma I_\sigma(u) = \rho \mu F(u)$, where $F(u) = \int_0^1 f(\tau u) u\, d\tau$ is the potential of $f$.
\end{enumerate}
Solutions of equation \eqref{18} are the critical points of the $C^1$-functional
\[
E(u) = I_\sigma(u) - F(u), \quad u \in W,
\]
where $I_\sigma$ is as in the last section. However, it seems difficult to verify the boundedness of \PS{} sequences of $E$ in general.

Let $\M$ be as in the last section. As noted in \cite[Section 2.2.1]{MR5043800}, the scaling induces a continuous projection
\[
W \setminus \set{0} \to \M, \quad u \mapsto \widetilde{u} := u_{t_u},
\]
where $t_u = I_\sigma(u)^{- 1/\sigma}$ (see \eqref{16}). This projection can be used to reduce equation \eqref{18} to the nonlinear eigenvalue problem
\begin{equation} \label{17}
\As(u) = \mu f(u)
\end{equation}
on $\M$. Using this constrained approach, the following theorem was proved in Perera \cite{Pe25}.

\begin{theorem}[{\cite[Theorem 3.8]{Pe25}}] \label{Theorem 4}
Assume $(A_1)$--$(A_7)$, $(A_{10})$, and $(f_1)$--$(f_4)$. Then equation \eqref{18} has an infinite sequence of nontrivial solutions
\begin{equation} \label{23}
\overline{u}_k = u^k_{\mu_k^{- 1/(\sigma - \rho)}},
\end{equation}
where $u^k \in \M$, $\mu_k > 0$, and $\mu_k \nearrow \infty$. Moreover,
\[
E(\overline{u}_k) = \left(1 - \frac{\sigma}{\rho}\right) \mu_k^{- \sigma/(\sigma - \rho)},
\]
in particular, $E(\overline{u}_k) < 0$ and $E(\overline{u}_k) \nearrow 0$ if $\rho < \sigma$, while $E(\overline{u}_k) > 0$ and $E(\overline{u}_k) \nearrow \infty$ if $\rho > \sigma$.
\end{theorem}

\section{Proofs of main theorems}

Throughout this section, for $u \in W$ we write
\[
K_i(u) = [u]_i^{p_i} = \int_{\R^N} \int_{\R^N} \frac{|u(x) - u(y)|^{p_i}}{|x - y|^{N + s_i p_i}}\, dx dy, \quad i = 1, 2, \qquad L(u) = \int_{\R^N} |x|^{-b}\, |u|^q\, dx,
\]
and for $u, v \in W$,
\[
\cE_i(u,v) = \int_{\R^N} \int_{\R^N} \frac{|u(x) - u(y)|^{p_i - 2}\, (u(x) - u(y))(v(x) - v(y))}{|x - y|^{N + s_i p_i}}\, dx dy, \quad i = 1, 2.
\]

\subsection{Scaling and the scaled operator}

Set
\[
u_t(x) = t^\beta\, u(tx), \quad x \in \R^N,\, t \ge 0.
\]
The changes of variables $x' = tx$, $y' = ty$ and \eqref{alpha-beta} give
\begin{equation} \label{scaling Ki}
K_i(u_t) = t^{\, p_i\, (\beta + s_i) - N}\, K_i(u) = t^\sigma K_i(u), \qquad [u_t]_i = t^{\sigma/p_i}\, [u]_i, \qquad \sigma = \alpha - N,
\end{equation}
and
\begin{equation} \label{scaling L}
L(u_t) = t^{\, \beta q + b - N}\, L(u) = t^\rho L(u), \qquad \rho = b - N + \beta q.
\end{equation}
Then the mapping
\begin{equation} \label{14}
W_\rad \times [0,\infty) \to W_\rad, \quad (u,t) \mapsto u_t
\end{equation}
is a scaling in $W_\rad$. Assumptions $(A_1)$--$(A_3)$ are clearly satisfied, while $(A_4)$ and $(A_5)$ follow from the estimate
\[
\norm{u_t} = t^{\sigma/p_1}\, [u]_1 + t^{\sigma/p_2}\, [u]_2 \le \max \set{t^{\sigma/p_1},t^{\sigma/p_2}} \norm{u}
\]
(see \eqref{12}). It only remains to verify continuity, which we do in the following proposition.

\begin{proposition}
The mapping in \eqref{14} is continuous.
\end{proposition}

\begin{proof}
Let $u_j \to u$ in $W_\rad$ and $t_j \to t$ in $[0,\infty)$. If $t = 0$, then
\[
\norm{(u_j)_{t_j}} \le \max \set{t_j^{\sigma/p_1},t_j^{\sigma/p_2}} \norm{u_j} \to 0 = \norm{u_0}.
\]
Let $t > 0$. Since
\[
\norm{(u_j)_{t_j} - u_t} \le \norm{(u_j - u)_{t_j}} + \norm{u_{t_j} - u_t} \le \max \set{t_j^{\sigma/p_1},t_j^{\sigma/p_2}} \norm{u_j - u} + \norm{u_{t_j} - u_t},
\]
it suffices to show that $[u_{t_j} - u_t]_i \to 0$ for $i = 1, 2$. First let $\varphi \in C^\infty_c(\R^N)$ and set
\[
\Phi_\tau(x,y) = \frac{\varphi_\tau(x) - \varphi_\tau(y)}{|x - y|^{(N + s_i p_i)/p_i}}, \quad \tau > 0.
\]
Then $\Phi_{t_j} \to \Phi_t$ pointwise on $\set{x \ne y}$, and $\norm[L^{p_i}(\R^{2N})]{\Phi_{t_j}} = t_j^{\sigma/p_i}\, [\varphi]_i \to t^{\sigma/p_i}\, [\varphi]_i = \norm[L^{p_i}(\R^{2N})]{\Phi_t}$ by \eqref{scaling Ki}, so $\Phi_{t_j} \to \Phi_t$ in $L^{p_i}(\R^{2N})$ by the Brezis-Lieb lemma \cite{MR699419}, i.e., $[\varphi_{t_j} - \varphi_t]_i \to 0$. For general $u \in W$ and $\eps > 0$, choose $\varphi \in C^\infty_c(\R^N)$ with $[u - \varphi]_i < \eps$, which is possible since $\dot{W}^{s_i,p_i}(\R^N)$ is the completion of $C^\infty_c(\R^N)$. Then
\[
[u_{t_j} - u_t]_i \le [(u - \varphi)_{t_j}]_i + [\varphi_{t_j} - \varphi_t]_i + [(\varphi - u)_t]_i \le \left(t_j^{\sigma/p_i} + t^{\sigma/p_i}\right) \eps + [\varphi_{t_j} - \varphi_t]_i,
\]
so $\limsup [u_{t_j} - u_t]_i \le 2 t^{\sigma/p_i}\, \eps$, and the conclusion follows since $\eps > 0$ is arbitrary.
\end{proof}

The scaled operator that corresponds to the left-hand side of equation \eqref{1} and its potential are given by
\[
\As(u) v = \cE_1(u,v) + \cE_2(u,v), \qquad I_\sigma(u) = \frac{1}{p_1}\, K_1(u) + \frac{1}{p_2}\, K_2(u)
\]
for $u, v \in W_\rad$. Assumptions $(A_6)$ and $(A_{10})$ are clearly satisfied. We verify assumption $(A_7)$ in the following proposition.

\begin{proposition} \label{Proposition S+}
If $u_j \wto u$ in $W$ and $\As(u_j)(u_j - u) \to 0$, then $u_j \to u$ in $W$.
\end{proposition}

\begin{proof}
Since $\As(u) \in W^\ast$ and $u_j - u \wto 0$, $\As(u)(u_j - u) \to 0$, so
\[
\sum_{i=1}^2 \big[\cE_i(u_j,u_j - u) - \cE_i(u,u_j - u)\big] \to 0.
\]
Write $J_p(r) = |r|^{p-2}\, r$, $U_j(x,y) = u_j(x) - u_j(y)$, $U(x,y) = u(x) - u(y)$, and $d\mu_i = |x - y|^{- (N + s_i p_i)}\, dx dy$. Each bracket equals $\iint \big(J_{p_i}(U_j) - J_{p_i}(U)\big)(U_j - U)\, d\mu_i$, which is nonnegative by the monotonicity of $J_{p_i}$, so each bracket tends to zero. If $p_i \ge 2$, the elementary inequality $(J_p(a) - J_p(b))(a - b) \ge 2^{2-p}\, |a - b|^p$ gives $K_i(u_j - u) = \iint |U_j - U|^{p_i}\, d\mu_i \to 0$. If $1 < p_i < 2$, the inequality $(J_p(a) - J_p(b))(a - b) \ge (p - 1)\, |a - b|^2\, (|a| + |b|)^{p-2}$ and the H\"{o}lder inequality give
\[
\iint |U_j - U|^{p_i}\, d\mu_i \le \left(\iint \frac{|U_j - U|^2}{(|U_j| + |U|)^{2 - p_i}}\, d\mu_i\right)^{p_i/2} \left(\iint (|U_j| + |U|)^{p_i}\, d\mu_i\right)^{(2 - p_i)/2} \to 0,
\]
since the second factor is bounded. Hence $\norm{u_j - u} \to 0$.
\end{proof}

\begin{lemma}[symmetric criticality] \label{Lemma symmetric criticality}
Let $T \in W^\ast$ satisfy $T(v \circ R) = T(v)$ for all $v \in W$ and $R \in \O(N)$. If $T(v) = 0$ for all $v \in W_\rad$, then $T = 0$. In particular, weak radial solutions of the equations considered in this paper are weak solutions.
\end{lemma}

\begin{proof}
The orthogonal group $\O(N)$ acts on $W$ by the linear isometries $v \mapsto v \circ R$. For $v \in W$, the radial average $\widetilde{v} = \int_{\O(N)} v \circ R\, dR$, taken with respect to the normalized Haar measure and converging as a Bochner integral in $W$, belongs to $W_\rad$, so $T(v) = \int_{\O(N)} T(v \circ R)\, dR = T(\widetilde{v}\thinspace) = 0$ (cf.\! \cite{MR547524}). The last statement follows since, for a radial $u$, the difference of the two sides of the weak formulation is a rotation-invariant element of $W^\ast$.
\end{proof}

\subsection{Embeddings into weighted Lebesgue spaces} \label{Embeddings}

Recall the sets $\T$, $\Tc$, and $\Te$ defined in the introduction. We will use the two families of Sobolev-invariant dilations
\[
v^{(i)}_\lambda(x) = \lambda^{(N - s_i p_i)/p_i}\, v(\lambda x), \quad \lambda > 0,\ i = 1, 2,
\]
which satisfy $[v^{(i)}_\lambda]_i = [v]_i$ and $\norm[L^{p_i^\ast}(\R^N)]{v^{(i)}_\lambda} = \norm[L^{p_i^\ast}(\R^N)]{v}$. A change of variables shows that
\begin{equation} \label{kappa}
[v^{(1)}_\lambda]_2 = \lambda^\kappa\, [v]_2, \qquad [v^{(2)}_\lambda]_1 = \lambda^{- \kappa}\, [v]_1,
\end{equation}
where
\[
\kappa := (s_2 - s_1) - N \left(\frac{1}{p_2} - \frac{1}{p_1}\right) = \frac{(p_1 - p_2)(\alpha - N)}{p_1 p_2} > 0.
\]

\begin{theorem} \label{Theorem embedding}
Let $q \ge 1$ and $b \ge 0$.
\begin{enumroman}
\item \label{Theorem embedding.i} If $(q,b) \in \Te$, then $W$ is continuously embedded in $L^q(\R^N;|x|^{-b})$.
\item \label{Theorem embedding.ii} If $(q,b) \in \Tc$ and $b > 0$, then $W$ is compactly embedded in $L^q(\R^N;|x|^{-b})$. If $N \ge 2$ and $p_1^\ast < q < p_2^\ast$, then $W_\rad$ is compactly embedded in $L^q(\R^N)$.
\item \label{Theorem embedding.iii} If $(q,b) \in [Q_2,Q_1] \cup [Q_3,Q_4]$, then the embedding $W_\rad \hookrightarrow L^q(\R^N;|x|^{-b})$ is not compact.
\end{enumroman}
\end{theorem}

We first record the elementary estimates on which the proof is based. For $0 < \eps < R$ write $B_\eps = \set{|x| < \eps}$ and $B_R^c = \set{|x| > R}$.

\begin{lemma} \label{Lemma splitting}
Let $q \ge 1$, $b \ge 0$, and $u \in W$.
\begin{enumroman}
\item If $q < p_2^\ast$ and $b < N\, (1 - q/p_2^\ast)$, then
\[
\int_{B_\eps} |x|^{-b}\, |u|^q\, dx \le C\, \eps^{N (1 - q/p_2^\ast) - b}\, [u]_2^q.
\]
\item If $q \le p_1^\ast$ and $b > N\, (1 - q/p_1^\ast)$, then
\[
\int_{B_R^c} |x|^{-b}\, |u|^q\, dx \le C\, R^{- (b - N (1 - q/p_1^\ast))}\, [u]_1^q.
\]
\item If $p_1^\ast \le q \le p_2^\ast$, then
\[
\int_{B_R^c} |x|^{-b}\, |u|^q\, dx \le R^{-b}\, \norm[L^{p_1^\ast}(\R^N)]{u}^{\theta q}\, \norm[L^{p_2^\ast}(\R^N)]{u}^{(1 - \theta)\, q} \le C\, R^{-b} \norm{u}^q, \qquad \frac{1}{q} = \frac{\theta}{p_1^\ast} + \frac{1 - \theta}{p_2^\ast}.
\]
\end{enumroman}
Here the constants $C$ depend only on $N$, $s_i$, $p_i$, $q$, and $b$.
\end{lemma}

\begin{proof}
$(i)$ By the H\"{o}lder inequality with exponents $p_2^\ast/(p_2^\ast - q)$ and $p_2^\ast/q$,
\[
\int_{B_\eps} |x|^{-b}\, |u|^q\, dx \le \left(\int_{B_\eps} |x|^{- b p_2^\ast/(p_2^\ast - q)}\, dx\right)^{1 - q/p_2^\ast} \norm[L^{p_2^\ast}(\R^N)]{u}^q,
\]
and the last integral is finite and equal to $C \eps^{N - b p_2^\ast/(p_2^\ast - q)}$ since $b p_2^\ast/(p_2^\ast - q) < N$. Now apply the fractional Sobolev inequality. $(ii)$ If $q < p_1^\ast$, argue as in $(i)$, the integral $\int_{B_R^c} |x|^{- b p_1^\ast/(p_1^\ast - q)}\, dx$ being finite and equal to $C R^{N - b p_1^\ast/(p_1^\ast - q)}$ since $b p_1^\ast/(p_1^\ast - q) > N$. If $q = p_1^\ast$, then $b > 0$ and the estimate follows from $|x|^{-b} \le R^{-b}$ on $B_R^c$. $(iii)$ follows from $|x|^{-b} \le R^{-b}$ on $B_R^c$, the H\"{o}lder interpolation inequality, and the fractional Sobolev inequalities.
\end{proof}

For radial functions, the tails at infinity are also uniformly small at the unweighted points $(q,0)$ with $p_1^\ast < q < p_2^\ast$. The following lemma does not require any pointwise decay estimate for radial functions, and hence no restriction of the form $s_i p_i > 1$.

\begin{lemma} \label{Lemma radial tail}
Let $N \ge 2$ and $p_1^\ast < q < p_2^\ast$. Then
\[
\sup_{u \in W_\rad,\, \norm{u} \le 1}\, \int_{B_R^c} |u|^q\, dx \to 0 \quad \text{as } R \to \infty.
\]
\end{lemma}

\begin{proof}
For $y \in \R^N$, write $B(y) = \set{|x - y| < 1}$ and
\[
[u]_{2,B(y)}^{p_2} = \int_{B(y)} \int_{B(y)} \frac{|u(x) - u(z)|^{p_2}}{|x - z|^{N + s_2 p_2}}\, dx dz.
\]
\emph{Step 1: local estimate.} By Jensen's inequality, $|u(x) - u_{B(y)}|^{p_2} \le \frac{1}{|B(y)|} \int_{B(y)} |u(x) - u(z)|^{p_2}\, dz \le C \int_{B(y)} |u(x) - u(z)|^{p_2}\, |x - z|^{- N - s_2 p_2}\, dz$ for $x \in B(y)$, where $u_{B(y)}$ is the mean of $u$ on $B(y)$, so $\norm[L^{p_2}(B(y))]{u - u_{B(y)}} \le C\, [u]_{2,B(y)}$. Combining this with the fractional Sobolev inequality in $W^{s_2,p_2}(B(y))$ (see, e.g., \cite[Theorem 6.7]{MR2944369} or \cite[Theorem 1.3]{DPSV26}) and the H\"{o}lder estimate $|u_{B(y)}|\, |B(y)|^{1/p_2^\ast} \le C \norm[L^{p_1^\ast}(B(y))]{u}$ gives
\begin{equation} \label{local Sobolev}
\norm[L^{p_2^\ast}(B(y))]{u} \le C \left([u]_{2,B(y)} + \norm[L^{p_1^\ast}(B(y))]{u}\right)
\end{equation}
with $C$ independent of $y$.

\emph{Step 2: decay of local norms.} Let $u \in W_\rad$ and $|y| = r \ge 2$. Since $N \ge 2$, there are rotations $R_1, \dots, R_n \in \O(N)$ with $n \ge c\, r^{N-1}$, $c = c(N) > 0$, such that the balls $B(R_k y)$ are pairwise disjoint. Since $u$ is radial, the quantities $\int_{B(R_k y)} |u|^{p_1^\ast}$ and $[u]_{2,B(R_k y)}^{p_2}$ do not depend on $k$, and since the sets $B(R_k y) \times B(R_k y)$ are pairwise disjoint,
\begin{equation} \label{radial local}
\int_{B(y)} |u|^{p_1^\ast}\, dx \le \frac{C}{r^{N-1}}\, \norm[L^{p_1^\ast}(\R^N)]{u}^{p_1^\ast}, \qquad [u]_{2,B(y)}^{p_2} \le \frac{C}{r^{N-1}}\, [u]_2^{p_2}.
\end{equation}

\emph{Step 3: summation.} Cover $B_R^c$, $R \ge 3$, by unit balls $B(y_k)$ with $|y_k| \ge R - 1$ and with at most $M = M(N)$ of them containing any given point. Then also $\sum_k \chi_{B(y_k) \times B(y_k)}(x,z) \le M$. Let $\theta \in (0,1)$ be given by $1/q = \theta/p_1^\ast + (1 - \theta)/p_2^\ast$, so that $\theta q/p_1^\ast + (1 - \theta)\, q/p_2^\ast = 1$. By the H\"{o}lder inequality on each ball and then for sums,
\begin{align*}
\int_{B_R^c} |u|^q\, dx & \le \sum_k \norm[L^{p_1^\ast}(B(y_k))]{u}^{\theta q}\, \norm[L^{p_2^\ast}(B(y_k))]{u}^{(1 - \theta)\, q}\\[5pt]
& \le \Big(\sum_k \norm[L^{p_1^\ast}(B(y_k))]{u}^{p_1^\ast}\Big)^{\theta q/p_1^\ast} \Big(\sum_k \norm[L^{p_2^\ast}(B(y_k))]{u}^{p_2^\ast}\Big)^{(1 - \theta)\, q/p_2^\ast}.
\end{align*}
The first factor is bounded by $(M \norm[L^{p_1^\ast}(\R^N)]{u}^{p_1^\ast})^{\theta q/p_1^\ast} \le C$ for $\norm{u} \le 1$. By \eqref{local Sobolev} and \eqref{radial local}, and since $p_2^\ast > p_2$ and $p_2^\ast > p_1^\ast$,
\begin{multline*}
\sum_k \norm[L^{p_2^\ast}(B(y_k))]{u}^{p_2^\ast} \le C \sup_k\, [u]_{2,B(y_k)}^{p_2^\ast - p_2} \sum_k [u]_{2,B(y_k)}^{p_2} + C \sup_k\, \norm[L^{p_1^\ast}(B(y_k))]{u}^{p_2^\ast - p_1^\ast} \sum_k \norm[L^{p_1^\ast}(B(y_k))]{u}^{p_1^\ast}\\[5pt]
\le C \left(R^{- (N-1)(p_2^\ast - p_2)/p_2} + R^{- (N-1)(p_2^\ast - p_1^\ast)/p_1^\ast}\right)
\end{multline*}
for $\norm{u} \le 1$. Since $\theta < 1$, the conclusion follows.
\end{proof}

\begin{proof}[Proof of Theorem \ref{Theorem embedding}]
\ref{Theorem embedding.i} First let $(q,b) \in \Tc$. Then $q < p_2^\ast$ and $b < N\, (1 - q/p_2^\ast)$, so the integral over $B_1$ is estimated by Lemma \ref{Lemma splitting} $(i)$. If $q \le p_1^\ast$, then $b > N\, (1 - q/p_1^\ast)$ and the integral over $B_1^c$ is estimated by Lemma \ref{Lemma splitting} $(ii)$; if $q > p_1^\ast$, it is estimated by Lemma \ref{Lemma splitting} $(iii)$. In all cases $\norm[L^q(\R^N;|x|^{-b})]{u} \le C \norm{u}$.

Next let $(q,b) \in [Q_2,Q_1]$. The map $(q,b) \mapsto (1/q,b/q)$ maps line segments in the half-plane $\set{q > 0}$ to line segments, so there is $\eta \in [0,1]$ with
\[
\frac{1}{q} = \frac{1 - \eta}{p_1} + \frac{\eta}{p_1^\ast}, \qquad \frac{b}{q} = (1 - \eta)\, s_1.
\]
By the H\"{o}lder inequality with exponents $p_1/((1 - \eta)\, q)$ and $p_1^\ast/(\eta q)$,
\[
\int_{\R^N} |x|^{-b}\, |u|^q\, dx \le \left(\int_{\R^N} \frac{|u|^{p_1}}{|x|^{s_1 p_1}}\, dx\right)^{(1 - \eta)\, q/p_1} \left(\int_{\R^N} |u|^{p_1^\ast}\, dx\right)^{\eta q/p_1^\ast} \le C\, [u]_1^q
\]
by the fractional Hardy and Sobolev inequalities. The case $(q,b) \in [Q_3,Q_4]$ is similar.

\ref{Theorem embedding.ii} Let $u_j \wto 0$ in $W$ (resp.\! $W_\rad$), so that $\seq{u_j}$ is bounded. Given $\eps > 0$ and $R > 1$, split
\[
\int_{\R^N} |x|^{-b}\, |u_j|^q\, dx = \int_{B_\eps} + \int_{B_R \setminus B_\eps} + \int_{B_R^c}.
\]
By Lemma \ref{Lemma splitting} $(i)$, the first integral is bounded by $C \eps^{N (1 - q/p_2^\ast) - b}$ uniformly in $j$. If $b > 0$, the third integral is bounded uniformly in $j$ by $C R^{- (b - N (1 - q/p_1^\ast))}$ when $q \le p_1^\ast$ and by $C R^{-b}$ when $q > p_1^\ast$, by Lemma \ref{Lemma splitting} $(ii)$ and $(iii)$; if $b = 0$, $N \ge 2$, and $u_j \in W_\rad$, it tends to zero as $R \to \infty$ uniformly in $j$ by Lemma \ref{Lemma radial tail}. Finally, $\seq{u_j}$ is bounded in $W^{s_2,p_2}(B_R)$ since $\norm[L^{p_2}(B_R)]{u_j} \le C \norm[L^{p_2^\ast}(\R^N)]{u_j}$, so $u_j \to 0$ in $L^q(B_R)$ for a renamed subsequence by the compactness of the embedding $W^{s_2,p_2}(B_R) \hookrightarrow L^q(B_R)$ for $q < p_2^\ast$ (see \cite[Corollary 7.2]{MR2944369} or \cite[Corollary 1.5]{DPSV26}), the limit being zero since $u_j \wto 0$; as $|x|^{-b} \le \eps^{-b}$ on $B_R \setminus B_\eps$, the second integral tends to zero. Hence $\int_{\R^N} |x|^{-b}\, |u_j|^q\, dx \to 0$ for a renamed subsequence, and since every subsequence of $\seq{u_j}$ has such a further subsequence, for the whole sequence.

\ref{Theorem embedding.iii} Let $(q,b) \in [Q_2,Q_1]$, so that $b = N\, (1 - q/p_1^\ast)$. Fix a nonzero radial $\varphi \in C^\infty_c(\R^N)$ and let $v_\lambda = \varphi^{(1)}_\lambda$ with $\lambda \searrow 0$. By \eqref{kappa}, $[v_\lambda]_1 = [\varphi]_1$ and $[v_\lambda]_2 = \lambda^\kappa\, [\varphi]_2 \to 0$, so $\seq{v_\lambda}$ is bounded in $W_\rad$, and $v_\lambda \wto 0$ since $v_\lambda \to 0$ uniformly. On the other hand, a change of variables gives
\[
\int_{\R^N} |x|^{-b}\, |v_\lambda|^q\, dx = \lambda^{q N/p_1^\ast + b - N} \int_{\R^N} |x|^{-b}\, |\varphi|^q\, dx = \int_{\R^N} |x|^{-b}\, |\varphi|^q\, dx > 0,
\]
so $\seq{v_\lambda}$ has no subsequence converging in $L^q(\R^N;|x|^{-b})$. If $(q,b) \in [Q_3,Q_4]$, argue in the same way with $v_\lambda = \varphi^{(2)}_\lambda$ and $\lambda \to \infty$, using $[v_\lambda]_1 = \lambda^{-\kappa}\, [\varphi]_1 \to 0$.
\end{proof}

\begin{remark} \label{Remark sharpness}
The region $\Te$ is optimal. If $b < N\, (1 - q/p_1^\ast)$, $b > N\, (1 - q/p_2^\ast)$, or $(q,b)$ lies on one of the half-open segments $[P,Q_2)$ and $[P,Q_3)$, then the embedding fails even for radial functions. Indeed, suitable functions behaving like $|x|^{-N/p_1^\ast}\, (\log |x|)^{-\gamma}$ near infinity, with $1/p_1 < \gamma \le 1/q$, belong to $W_\rad$ but not to $L^q(\R^N;|x|^{-b})$ when $b = N\, (1 - q/p_1^\ast)$ and $q < p_1$, and similarly with $|x|^{-N/p_2^\ast}\, |\log |x||^{-\gamma}$ near the origin on $[P,Q_3)$. Moreover, since the points of $\Tc$ with $b = 0$ are translation-invariant, the restriction to radial functions is essential for the compactness of the embedding on the bottom edge $(Q_1,Q_4)$ of $\T$.
\end{remark}

\begin{remark}
In contrast with the corresponding result for the double phase operator in \cite{PeWi}, the continuous embeddings in Theorem \ref{Theorem embedding} hold in the whole space $W$, and the compact ones hold in $W$ off the bottom edge of $\T$: the two fractional Sobolev inequalities hold in $\dot{W}^{s_i,p_i}(\R^N)$ without symmetry, and the weight $|x|^{-b}$ controls the region near infinity when $b > 0$.
\end{remark}

\subsection{Weak limits of Palais-Smale sequences}

For the critical problem \eqref{BN} we will need the following lemma.

\begin{lemma} \label{Lemma weak limit}
Assume \eqref{BN hypotheses} and let $\seq{u_j}$ be a bounded sequence in $W_\rad$ such that $E_\mu'(u_j) \to 0$. Then a renamed subsequence converges weakly in $W_\rad$ and a.e.\! in $\R^N$ to a critical point $u$ of $E_\mu$.
\end{lemma}

\begin{proof}
Since $W_\rad$ is reflexive, a renamed subsequence converges weakly to some $u \in W_\rad$. By the compactness of the embeddings $W^{s_2,p_2}(B_R) \hookrightarrow L^1(B_R)$ and a diagonal argument, a further renamed subsequence converges to $u$ a.e.\! in $\R^N$. Let $v \in W_\rad$. For $i = 1, 2$, the functions
\[
J_{p_i}(u_j(x) - u_j(y))\, |x - y|^{- (N + s_i p_i)/p_i'}
\]
are bounded in $L^{p_i'}(\R^{2N})$ and converge a.e.\! to $J_{p_i}(u(x) - u(y))\, |x - y|^{- (N + s_i p_i)/p_i'}$, and hence weakly in $L^{p_i'}(\R^{2N})$. Since $(v(x) - v(y))\, |x - y|^{- (N + s_i p_i)/p_i} \in L^{p_i}(\R^{2N})$, it follows that $\cE_i(u_j,v) \to \cE_i(u,v)$. Similarly, for $r \in \set{q,p_1^\ast,p_2^\ast}$, the functions $|u_j|^{r-2}\, u_j$ are bounded in $L^{r'}(\R^N)$ and converge a.e., and hence weakly in $L^{r'}(\R^N)$, to $|u|^{r-2}\, u$, so $\int_{\R^N} |u_j|^{r-2}\, u_j v\, dx \to \int_{\R^N} |u|^{r-2}\, uv\, dx$. Hence $E_\mu'(u)\, v = \lim E_\mu'(u_j)\, v = 0$.
\end{proof}

\subsection{Poho\v{z}aev identity} \label{Subsection Pohozaev proof}

In this subsection we prove Theorem \ref{Theorem Pohozaev}. The key analytical ingredient is the following estimate for the dilation modulus of a fractional Sobolev function in a Hardy-weighted norm. Note that the displacement $|tx - x| = |t - 1|\, |x|$ grows linearly in $|x|$, and that the homogeneous weight $|x|^{-sp}$ exactly compensates this growth, so that the full rate $|t - 1|^s$ is obtained without any relation between $s$, $p$, and $N$.

\begin{lemma} \label{Lemma modulus}
Let $s \in (0,1)$ and $p \in [1,\infty)$. Then there exists a constant $C = C(N,s,p) > 0$ such that for all $u \in \dot{W}^{s,p}(\R^N)$ and all $t > 0$ with $|t - 1| \le 1/12$,
\[
\int_{\R^N} \frac{|u(tx) - u(x)|^p}{|x|^{sp}}\, dx \le C\, |t - 1|^{sp}\, [u]_{\dot{W}^{s,p}(\R^N)}^p.
\]
\end{lemma}

\begin{proof}
Let $\eps = |t - 1| \le 1/12$. For $k \in \Z$ set
\[
A_k = \set{2^k \le |x| < 2^{k+1}}, \qquad \widetilde{A}_k = \set{2^{k-2} < |x| < 2^{k+3}}, \qquad \delta_k = 4 \eps 2^k,
\]
and denote by $[u]_{\widetilde{A}_k}^p$ the Gagliardo integral of $u$ over $\widetilde{A}_k \times \widetilde{A}_k$. For $x \in A_k$ we have $|tx - x| = \eps |x| < \delta_k/2$. For any $\xi \in B_{\delta_k}(x)$,
\[
|u(tx) - u(x)|^p \le 2^{p-1} \left(|u(tx) - u(\xi)|^p + |u(\xi) - u(x)|^p\right),
\]
and averaging over $\xi \in B_{\delta_k}(x)$ and integrating over $x \in A_k$ gives
\begin{multline*}
\int_{A_k} |u(tx) - u(x)|^p\, dx\\[5pt]
\le \frac{C}{\delta_k^N} \left[\iint_{\substack{x \in A_k\\ |\xi - x| < \delta_k}} |u(\xi) - u(x)|^p\, d\xi\, dx + \iint_{\substack{x \in A_k\\ |\xi - x| < \delta_k}} |u(tx) - u(\xi)|^p\, d\xi\, dx\right].
\end{multline*}
In the first integral, $\delta_k^{-N} \le \delta_k^{sp}\, |x - \xi|^{- N - sp}$, and all pairs $(x,\xi)$ lie in $\widetilde{A}_k \times \widetilde{A}_k$ since $|\xi| \ge 2^k (1 - 4 \eps) > 2^{k-2}$ and $|\xi| < 2^{k+1} + \delta_k < 2^{k+2}$. So the first term is at most $C \delta_k^{sp}\, [u]_{\widetilde{A}_k}^p$. In the second integral we substitute $y = tx$, so that $dx = t^{-N}\, dy \le C\, dy$, $y \in \widetilde{A}_k$, and $|\xi - y| \le |\xi - x| + |x - y| < 2 \delta_k$, so as before the second term is at most $C \delta_k^{sp}\, [u]_{\widetilde{A}_k}^p$. Multiplying by $\sup_{A_k} |x|^{-sp} = 2^{- spk}$ and using $2^{- spk}\, \delta_k^{sp} = (4 \eps)^{sp}$ gives
\[
\int_{A_k} \frac{|u(tx) - u(x)|^p}{|x|^{sp}}\, dx \le C \eps^{sp}\, [u]_{\widetilde{A}_k}^p.
\]
Summing over $k \in \Z$ and noting that each pair $(x,y)$ belongs to $\widetilde{A}_k \times \widetilde{A}_k$ for at most five values of $k$ completes the proof.
\end{proof}

At a general point of $\H$, interpolation converts Lemma \ref{Lemma modulus} into the rate $b/q$.

\begin{lemma} \label{Lemma hull modulus}
Let $(q,b) \in \H$. Then there exists a constant $C > 0$ such that for all $u \in W$ and all $t > 0$ with $|t - 1| \le 1/12$,
\[
\left(\int_{\R^N} |x|^{-b}\, |u(tx) - u(x)|^q\, dx\right)^{1/q} \le C\, |t - 1|^{b/q} \norm{u}.
\]
\end{lemma}

\begin{proof}
Write $Q_k = (q_k,b_k)$, $k = 1, \dots, 4$, so that $b_1 = b_4 = 0$, $b_2/q_2 = s_1$, and $b_3/q_3 = s_2$. Since the map $\Pi : (q,b) \mapsto (1/q,b/q)$ is a homeomorphism of the half-plane $\set{q > 0}$ onto itself that maps line segments to line segments, $\Pi(\H)$ is the convex hull of $\Pi(Q_1), \dots, \Pi(Q_4)$. Hence there are $\eta_k \ge 0$ with $\sum_k \eta_k = 1$ and
\begin{equation} \label{weights}
\frac{1}{q} = \sum_{k=1}^4 \frac{\eta_k}{q_k}, \qquad \frac{b}{q} = \eta_2\, s_1 + \eta_3\, s_2.
\end{equation}
By the H\"{o}lder inequality with exponents $q_k/(q \eta_k)$ over the indices with $\eta_k > 0$,
\[
\int_{\R^N} |x|^{-b}\, |w|^q\, dx = \int_{\R^N} \prod_k \left(|x|^{- b_k/q_k}\, |w|\right)^{q \eta_k} dx \le \prod_k \norm[L^{q_k}(\R^N)]{|x|^{- b_k/q_k}\, w}^{q \eta_k},
\]
the weights matching by \eqref{weights}. Apply this to $w = u(t\, \cdot) - u$. The factors with $k = 1, 4$ are bounded by $\norm[L^{p_i^\ast}(\R^N)]{u(t\, \cdot)} + \norm[L^{p_i^\ast}(\R^N)]{u} \le C\, [u]_i$ by the fractional Sobolev inequalities, and those with $k = 2, 3$ by $C\, |t - 1|^{s_i}\, [u]_i$ by Lemma \ref{Lemma modulus}. The product carries the rate $|t - 1|^{q (\eta_2 s_1 + \eta_3 s_2)} = |t - 1|^b$.
\end{proof}

\begin{corollary} \label{Corollary scaled modulus}
Let $(q,b) \in \H$ and $u \in W$. Then, for $|t - 1| \le 1/12$,
\[
\omega(t) := \left(\int_{\R^N} |x|^{-b}\, |u_t - u|^q\, dx\right)^{1/q} \le C\, |t - 1|^{\min \set{b/q,\, 1}} \norm{u}.
\]
\end{corollary}

\begin{proof}
Write $u_t - u = t^\beta\, (u(t\, \cdot) - u) + (t^\beta - 1)\, u$. The first term is estimated by Lemma \ref{Lemma hull modulus}, and the second by $|t^\beta - 1| \le C\, |t - 1|$ and the continuous embedding $W \hookrightarrow L^q(\R^N;|x|^{-b})$.
\end{proof}

\begin{proof}[Proof of Theorem \ref{Theorem Pohozaev}]
If $\lambda \le 0$, then testing the equation with $u$ gives $K_1(u) + K_2(u) = \lambda L(u) \le 0$, so $u = 0$ and \eqref{Pohozaev intro} is trivial. So let $\lambda > 0$.

\emph{Step 1.} For $t > 0$, the function $\varphi_t := u_t - u$ belongs to $W_\rad$, so it is an admissible test function:
\begin{equation} \label{tested}
\cE_1(u,\varphi_t) + \cE_2(u,\varphi_t) = \lambda \int_{\R^N} |x|^{-b}\, |u|^{q-2}\, u\, \varphi_t\, dx.
\end{equation}

\emph{Step 2: an exact expansion.} For $p > 1$ and $a, a' \in \R$, the convexity of $r \mapsto |r|^p$ gives
\[
0 \le R_p(a,a') := |a'|^p - |a|^p - p\, J_p(a)(a' - a) \le C_p \left(|a|^p + |a'|^p\right).
\]
Taking $a = u(x) - u(y)$ and $a' = u_t(x) - u_t(y)$ and integrating against $|x - y|^{- (N + s_i p_i)}\, dx dy$, all the resulting integrals being finite by \eqref{scaling Ki}, we obtain
\[
\cE_i(u,\varphi_t) = \frac{t^\sigma - 1}{p_i}\, K_i(u) - \frac{1}{p_i} \iint R_{p_i}\, \frac{dx dy}{|x - y|^{N + s_i p_i}}.
\]
Similarly, with $G(r) = |r|^q/q$ and
\[
0 \le S(a,a') := G(a') - G(a) - G'(a)(a' - a) \le C_q \left(|a|^q + |a'|^q\right),
\]
we obtain from \eqref{scaling L}
\[
\int_{\R^N} |x|^{-b}\, |u|^{q-2}\, u\, \varphi_t\, dx = \frac{t^\rho - 1}{q}\, L(u) - \mathcal{Q}_t, \qquad \mathcal{Q}_t := \int_{\R^N} |x|^{-b}\, S(u,u_t)\, dx \ge 0.
\]
Set
\[
\mathcal{P}_t := \sum_{i=1}^2 \frac{1}{p_i} \iint R_{p_i}\, \frac{dx dy}{|x - y|^{N + s_i p_i}} \ge 0, \qquad \Phi(t) := t^\sigma I_\sigma(u) - \frac{\lambda\, t^\rho}{q}\, L(u).
\]
Then \eqref{tested} can be rewritten as the exact identity
\begin{equation} \label{star}
\Phi(t) - \Phi(1) = \mathcal{P}_t - \lambda\, \mathcal{Q}_t \quad \forall t > 0.
\end{equation}

\emph{Step 3: the potential defect is superlinear.} Let $|t - 1| \le 1/12$. If $q \ge 2$, then $S(a,a') \le C_q\, (|a| + |a'|)^{q-2}\, (a' - a)^2$, so by the H\"{o}lder inequality with exponents $q/(q - 2)$ and $q/2$,
\[
\mathcal{Q}_t \le C \left(\int_{\R^N} |x|^{-b}\, (|u| + |u_t|)^q\, dx\right)^{(q-2)/q} \omega(t)^2 \le C\, \omega(t)^2,
\]
since $L(u_t) = t^\rho L(u)$ is bounded for $t$ near $1$. If $1 < q < 2$, then $S(a,a') \le C_q\, |a' - a|^q$ and $\mathcal{Q}_t \le C\, \omega(t)^q$. In both cases Corollary \ref{Corollary scaled modulus} gives
\[
\mathcal{Q}_t \le C\, |t - 1|^\varkappa, \qquad \varkappa := \min \set{b/q,\, 1} \cdot \min \set{q,2},
\]
and an elementary check shows that $\varkappa > 1$ if and only if $2b > \max \set{q,2}$, i.e., \eqref{C} holds.

\emph{Step 4: conclusion.} The function $\Phi$ is smooth on $(0,\infty)$ with
\[
\Phi'(1) = \sigma I_\sigma(u) - \frac{\rho \lambda}{q}\, L(u) =: \delta.
\]
By \eqref{star}, $\mathcal{P}_t \ge 0$, and Step 3,
\[
\Phi(t) - \Phi(1) \ge - \lambda\, \mathcal{Q}_t \ge - C\, |t - 1|^\varkappa \quad \text{for } |t - 1| \le 1/12.
\]
If $\delta \ne 0$, choose $t_n \to 1$ with $\delta\, (t_n - 1) < 0$. Then $\Phi(t_n) - \Phi(1) = \delta\, (t_n - 1) + \O((t_n - 1)^2) \le \frac{1}{2}\, \delta\, (t_n - 1)$ for large $n$, and the last display gives $\frac{1}{2}\, |\delta| \le C\, |t_n - 1|^{\varkappa - 1} \to 0$, a contradiction. Hence $\delta = 0$, which is \eqref{Pohozaev intro}.
\end{proof}

\subsection{Scaled eigenvalue problem and the proof of Theorem \ref{Theorem 1}} \label{Subsection proof Theorem 1}

Equation \eqref{15} with $p_2 < q < p_1$ and $b = \alpha - \beta q$ is a scaled eigenvalue problem with $\sigma = \alpha - N$. The scaled operator that corresponds to the right-hand side of equation \eqref{15} and its potential are given by
\[
\Bs(u) v = \int_{\R^N} |x|^{-b}\, |u|^{q-2}\, uv\, dx, \quad u, v \in W_\rad
\]
and
\[
J_\sigma(u) = \int_{\R^N} \frac{1}{q}\, |x|^{-b}\, |u|^q\, dx, \quad u \in W_\rad,
\]
respectively. Assumption $(A_8)$ is clearly satisfied, while $(A_9)$ follows from the compactness of the embedding $W_\rad \hookrightarrow L^q(\R^N;|x|^{-b})$ (see Theorem \ref{Theorem embedding} \ref{Theorem embedding.ii}), since $(q,b) \in \Tc$ and $b > 0$. Thus, $(A_1)$--$(A_{10})$ hold, but we do not assume \eqref{C} here, so we do not have $(A_{11})$ at our disposal.

Let $\M$, $\widetilde{\Psi}$, and $\lambda_k = \lambda_k(q)$ be as in Section \ref{Scaled eigenvalue problems}. Assumption $(A_{11})$ is used in \cite{MR5043800} only to identify the critical values of $\widetilde{\Psi}$ with eigenvalues (see \cite[Proposition 2.5]{MR5043800}). In particular, $\widetilde{\Psi}$ satisfies the \PS{} condition (see \cite[Proposition 2.6]{MR5043800}), and the proof of \cite[Theorem 2.10]{MR5043800} shows that $\lambda_k \nearrow \infty$ is a sequence of critical values of $\widetilde{\Psi}$. Let $u^k \in \M$ be a critical point of $\widetilde{\Psi}$ with $\widetilde{\Psi}(u^k) = \lambda_k$. Then
\[
\mu_k \As(u^k) = \lambda_k^2\, \Bs(u^k)
\]
for some $\mu_k \in \R$ (see \cite[(2.13)]{MR5043800}). Testing with $u^k$ gives $\mu_k \As(u^k)\, u^k = \lambda_k^2\, \Bs(u^k)\, u^k > 0$, so $\mu_k > 0$ by $(A_6)$ and $(A_8)$, and hence $u^k$ is a radial eigenfunction associated with the eigenvalue
\[
\widehat{\lambda}_k := \frac{\lambda_k^2}{\mu_k} = \frac{\As(u^k)\, u^k}{\Bs(u^k)\, u^k} = \frac{K_1(u^k) + K_2(u^k)}{L(u^k)}.
\]
Since $u^k \in \M$,
\[
\frac{1}{p_1}\, K_1(u^k) + \frac{1}{p_2}\, K_2(u^k) = I_\sigma(u^k) = 1, \qquad \frac{1}{q}\, L(u^k) = J_\sigma(u^k) = \frac{1}{\lambda_k},
\]
and since $p_2 < p_1$, the first equation gives $p_2 \le K_1(u^k) + K_2(u^k) \le p_1$. Hence
\[
\frac{p_2}{q}\, \lambda_k \le \widehat{\lambda}_k \le \frac{p_1}{q}\, \lambda_k,
\]
and therefore $\widehat{\lambda}_k \to \infty$. This proves Theorem \ref{Theorem 1}.

If, in addition, $2b > \max \set{q,2}$, then $(A_{11})$ holds by Theorem \ref{Theorem Pohozaev} since $\rho = \sigma$ on the scaling line, so $I_\sigma(u^k) = \widehat{\lambda}_k J_\sigma(u^k)$, i.e., $\widehat{\lambda}_k = \lambda_k$. In this case Theorem \ref{Theorem 3} also applies.

Finally, we show that $\widehat{\lambda}_1(q)=\lambda_1(q)$ without assuming \eqref{C}. Since every symmetric
two-point subset of $\M$ has cohomological index one, the definition of $\lambda_1(q)$ gives
\[
\lambda_1(q) = \inf_{u\in \M}\widetilde{\Psi}(u).
\]
Let $u_1\in \M$ be a critical point of $\widetilde{\Psi}$ at level $\lambda_1(q)$, as obtained above. Then $u_1$ is a global minimizer of $\widetilde{\Psi}$ on $\M$.

For every $u\in W_\rad\setminus\{0\}$, the scaling $u_t$ with $t=I_\sigma(u)^{-1/\sigma}$ belongs to $\M$. Since $\rho=\sigma$, we have
\[
\frac{I_\sigma(u_t)}{J_\sigma(u_t)} = \frac{I_\sigma(u)}{J_\sigma(u)}.
\]
Accordingly,
\[
\lambda_1(q) = \inf_{u\in W_\rad\setminus\{0\}} \frac{I_\sigma(u)}{J_\sigma(u)},
\]
and this infimum is attained at $u_1$.  Since $u_1$ minimizes $I_\sigma/J_\sigma$, we have
\[
0=\left(\frac{I_\sigma}{J_\sigma}\right)'(u_1)v =\frac{1}{J_\sigma(u_1)} \left(A_\sigma(u_1)v-\lambda_1(q)B_\sigma(u_1)v\right)
\]
for all $v\in W_\rad$. Hence $A_\sigma(u_1)=\lambda_1(q)B_\sigma(u_1)$, which proves that $\widehat{\lambda}_1(q)=\lambda_1(q)$. 

Moreover, since
\[
I_\sigma(|u_1|)\le I_\sigma(u_1)\quad\mbox{ and }\quad J_\sigma(|u_1|)=J_\sigma(u_1),
\]
the function $|u_1|$ also minimizes $I_\sigma/J_\sigma$, and hence is also an eigenfunction associated with $\lambda_1(q)$ by the argument above. So, replacing $u_1$ with $|u_1|$, we may assume that $u_1 \ge 0$.

Finally, by the definition of $\lambda_1(q)$ and by testing the eigenvalue equation with $u_1$, we obtain that
\[
\frac{K_1(u_1)}{p_1} + \frac{K_2(u_1)}{p_2} = \frac{\lambda_1(q)}{q}L(u_1) = \frac{K_1(u_1)+K_2(u_1)}{q}.
\]
Therefore,
\[
\left(\frac1q-\frac1{p_1}\right)K_1(u_1) = \left(\frac1{p_2}-\frac1q\right)K_2(u_1),
\]
and hence
\[
\frac{K_1(u_1)}{K_2(u_1)} = \frac{p_1(q-p_2)}{p_2(p_1-q)},
\]
as claimed.

\subsection{Proof of Theorem \ref{Theorem 2}}

By Theorem \ref{Theorem embedding} \ref{Theorem embedding.ii}, the embedding $W_\rad \hookrightarrow L^q(\R^N;|x|^{-b})$ is compact.

\emph{Case $(i)$: $q > p_1$.} In this case we do not use the scaling. We apply the symmetric mountain pass theorem of Ambrosetti and Rabinowitz \cite{MR0370183} (see also \cite[Theorem 9.12]{MR845785}) to the even functional $E$ on $W_\rad$, with $E(0) = 0$. Since
\[
E(u) - \frac{1}{q}\, E'(u)\, u = \left(\frac{1}{p_1} - \frac{1}{q}\right) K_1(u) + \left(\frac{1}{p_2} - \frac{1}{q}\right) K_2(u)
\]
and $p_2 < p_1 < q$, for a \PS{c} sequence $\seq{u_j}$ we have
\[
\left(\frac{1}{p_1} - \frac{1}{q}\right) \left([u_j]_1^{p_1} + [u_j]_2^{p_2}\right) \le c + \o(1) + \o(1) \norm{u_j},
\]
so $\seq{u_j}$ is bounded. For a renamed subsequence, $u_j \wto u$ in $W_\rad$, and hence $u_j \to u$ in $L^q(\R^N;|x|^{-b})$ by the compactness of the embedding. Since $\seq{|u_j|^{q-2}\, u_j}$ is bounded in $L^{q'}(\R^N;|x|^{-b})$, it follows that
\[
\As(u_j)(u_j - u) = E'(u_j)(u_j - u) + \lambda \int_{\R^N} |x|^{-b}\, |u_j|^{q-2}\, u_j\, (u_j - u)\, dx \to 0,
\]
and hence $u_j \to u$ by Proposition \ref{Proposition S+}. Thus $E$ satisfies the \PS{c} condition for all $c \in \R$. For $\norm{u} \le 1$ we have $[u]_i \le 1$, so
\[
E(u) \ge \frac{1}{p_1} \left([u]_1^{p_1} + [u]_2^{p_1}\right) - \frac{\lambda}{q}\, L(u) \ge \frac{1}{2^{p_1 - 1}\, p_1} \norm{u}^{p_1} - C \norm{u}^q,
\]
and since $q > p_1$, there exist $r, a > 0$ with $E \ge a$ on $\set{\norm{u} = r}$. On the other hand, on every finite-dimensional subspace $V$ of $W_\rad$ all norms are equivalent and
\[
E(u) \le \frac{1}{p_1} \norm{u}^{p_1} + \frac{1}{p_2} \norm{u}^{p_2} - \frac{\lambda}{q}\, L(u) \to - \infty \quad \text{as } \norm{u} \to \infty,\ u \in V,
\]
since $q > p_1 > p_2$ and $L^{1/q}$ is a norm on $V$. The symmetric mountain pass theorem gives an unbounded sequence of critical values $E(u_k) \ge a > 0$, and the $u_k$ are weak solutions by Lemma \ref{Lemma symmetric criticality}.

\emph{Case $(ii)$: $2b > \max \set{q,2}$.} We apply Theorem \ref{Theorem 4} to equation \eqref{15}. The scaled operator that corresponds to the right-hand side of equation \eqref{15} and its potential are given by
\[
f(u) v = \int_{\R^N} |x|^{-b}\, |u|^{q-2}\, uv\, dx, \quad u, v \in W_\rad
\]
and
\[
F(u) = \int_{\R^N} \frac{1}{q}\, |x|^{-b}\, |u|^q\, dx, \quad u \in W_\rad,
\]
respectively. Assumption $(f_1)$ holds with $\rho = b - N + \beta q > \sigma$ since $b > \alpha - \beta q$, $(f_2)$ is clearly satisfied, and $(f_3)$ follows from the compactness of the embedding $W_\rad \hookrightarrow L^q(\R^N;|x|^{-b})$. Assumption $(f_4)$ follows from Theorem \ref{Theorem Pohozaev}. Theorem \ref{Theorem 2} now follows from Theorem \ref{Theorem 4}.

\begin{proof}[Proof of Corollary \ref{Corollary pure power}]
The point $(q,0)$ with $p_1^\ast < q < p_2^\ast$ lies in $\H$, satisfies $N\, (1 - q/p_1^\ast) < 0 < N\, (1 - q/p_2^\ast)$ and $0 > \alpha - \beta q$, and $q > p_1$, so case $(i)$ of Theorem \ref{Theorem 2} applies.
\end{proof}

\subsection{Proof of Theorem \ref{Theorem nonexistence}}

Let $u \in W_\rad$ be a weak radial solution. If $\lambda \le 0$, then $u = 0$ as in the proof of Theorem \ref{Theorem Pohozaev}, so let $\lambda > 0$. Testing the equation with $u$ gives the Nehari identity
\begin{equation} \label{Nehari}
K_1(u) + K_2(u) = \lambda L(u).
\end{equation}
Let $(q,b) \in [Q_2,Q_1]$, so that $b = N\, (1 - q/p_1^\ast)$. By \eqref{alpha-beta} and $N/p_1^\ast = N/p_1 - s_1$,
\[
\rho = b - N + \beta q = q \left(\beta - \frac{N}{p_1^\ast}\right) = q \left(\frac{\alpha}{p_1} - s_1 - \frac{N}{p_1} + s_1\right) = \frac{\sigma q}{p_1}.
\]
By Theorem \ref{Theorem Pohozaev} and \eqref{Nehari},
\[
\sigma \left(\frac{K_1(u)}{p_1} + \frac{K_2(u)}{p_2}\right) = \frac{\rho \lambda}{q}\, L(u) = \frac{\sigma}{p_1} \left(K_1(u) + K_2(u)\right),
\]
so $\sigma\, (1/p_2 - 1/p_1)\, K_2(u) = 0$. Since $\sigma > 0$ and $p_2 < p_1$, $K_2(u) = 0$, and hence $u = 0$. If $(q,b) \in [Q_3,Q_4]$, then similarly $\rho = \sigma q/p_2$, and we obtain $\sigma\, (1/p_1 - 1/p_2)\, K_1(u) = 0$, so again $u = 0$. \qed

\subsection{Brezis-Nirenberg type problems: proof of Theorem \ref{Theorem BN}}

Throughout this subsection we assume \eqref{12} and \eqref{BN hypotheses}, and write
\[
X(u) = K_1(u), \qquad Y(u) = K_2(u), \qquad Z_i(u) = \int_{\R^N} |u|^{p_i^\ast}\, dx, \quad i = 1, 2.
\]
By the definition of $S_i$,
\begin{equation} \label{Sobolev Si}
Z_1(u) \le S_1^{- p_1^\ast/p_1}\, X(u)^{p_1^\ast/p_1}, \qquad Z_2(u) \le S_2^{- p_2^\ast/p_2}\, Y(u)^{p_2^\ast/p_2} \quad \forall u \in W.
\end{equation}
Note that
\begin{equation} \label{exponent order}
p_2 < p_1 < p_1^\ast < q < p_2^\ast, \qquad \frac{1}{p_i} - \frac{1}{p_i^\ast} = \frac{s_i}{N}.
\end{equation}

\subsubsection{Boundedness of Palais-Smale sequences} \label{PS:boundedness}

\begin{lemma} \label{Lemma bounded BN}
Every \PS{c} sequence of $E_\mu$ is bounded, and every critical point $u$ of $E_\mu$ satisfies $E_\mu(u) \ge 0$.
\end{lemma}

\begin{proof}
By \eqref{exponent order},
\begin{multline*}
E_\mu(u) - \frac{1}{p_1^\ast}\, E_\mu'(u)\, u = \left(\frac{1}{p_1} - \frac{1}{p_1^\ast}\right) X(u) + \left(\frac{1}{p_2} - \frac{1}{p_1^\ast}\right) Y(u) + \mu \left(\frac{1}{p_1^\ast} - \frac{1}{q}\right) \int_{\R^N} |u|^q\, dx\\[5pt]
+ \mu_2 \left(\frac{1}{p_1^\ast} - \frac{1}{p_2^\ast}\right) Z_2(u) \ge \frac{s_1}{N} \left(X(u) + Y(u)\right),
\end{multline*}
all the coefficients being positive. For a \PS{c} sequence $\seq{u_j}$, the left-hand side is $c + \o(1) + \o(1) \norm{u_j}$, which implies that $\seq{u_j}$ is bounded. For a critical point $u$, the left-hand side is $E_\mu(u)$.
\end{proof}

\begin{remark}
Lemma \ref{Lemma bounded BN} replaces the use of the Poho\v{z}aev identity in the corresponding part of \cite{PeWi}, where it was needed to show that $E(u) \ge 0$ for all critical points.
\end{remark}

\subsubsection{Threshold for compactness} \label{PS:threshold}

\begin{lemma} \label{Lemma threshold}
Let $X, Y, Z, W \ge 0$ satisfy $(X,Y) \ne (0,0)$ and
\begin{equation} \label{XYZW}
X + Y = Z + W, \qquad Z \le \mu_1 S_1^{- p_1^\ast/p_1}\, X^{p_1^\ast/p_1}, \qquad W \le \mu_2\, S_2^{- p_2^\ast/p_2}\, Y^{p_2^\ast/p_2}.
\end{equation}
Then
\begin{equation} \label{ineq}
I := \frac{1}{p_1}\, X + \frac{1}{p_2}\, Y - \frac{1}{p_1^\ast}\, Z - \frac{1}{p_2^\ast}\, W \ge \beta^\ast.
\end{equation}
This inequality is sharp, with equality for $(X,Y,Z,W) = (X_1,0,X_1,0)$ if $\beta^\ast = \beta_1$ and for $(0,Y_2,0,Y_2)$ if $\beta^\ast = \beta_2$, where
\[
X_1 = S_1^{N/s_1 p_1}\, \mu_1^{- (N - s_1 p_1)/s_1 p_1}, \qquad Y_2 = S_2^{N/s_2 p_2}\, \mu_2^{- (N - s_2 p_2)/s_2 p_2}.
\]
\end{lemma}

\begin{proof}
Since $X + Y = Z + W$, either $Z \ge X$ or $W \ge Y$.

\emph{Case 1: $Z \ge X$.} If $X = 0$, then $Z = 0$ by \eqref{XYZW}, so $W = Y > 0$ and we are in Case 2. So let $X > 0$. Then $Z > 0$, so $\mu_1 > 0$, and $X \le Z \le \mu_1 S_1^{- p_1^\ast/p_1}\, X^{p_1^\ast/p_1}$ gives $X \ge X_1$ since $1 - p_1/p_1^\ast = s_1 p_1/N$. Moreover, $W = X + Y - Z \le Y$ and
\[
I = \left(\frac{1}{p_1} - \frac{1}{p_1^\ast}\right) X + \left(\frac{1}{p_2} - \frac{1}{p_1^\ast}\right) Y + \left(\frac{1}{p_1^\ast} - \frac{1}{p_2^\ast}\right) W \ge \frac{s_1}{N}\, X \ge \frac{s_1}{N}\, X_1 = \beta_1,
\]
since the last two coefficients are positive by \eqref{exponent order}.

\emph{Case 2: $W \ge Y$.} If $Y = 0$, then $W = 0$, so $Z = X > 0$ and we are in Case 1. So let $Y > 0$. Then $W > 0$, so $\mu_2 > 0$, and as above $Y \ge Y_2$. Moreover, $Z = X + Y - W \le X$ and
\[
I = \left(\frac{1}{p_1} - \frac{1}{p_2^\ast}\right) X + \left(\frac{1}{p_2} - \frac{1}{p_2^\ast}\right) Y - \left(\frac{1}{p_1^\ast} - \frac{1}{p_2^\ast}\right) Z \ge \left(\frac{1}{p_1} - \frac{1}{p_1^\ast}\right) X + \frac{s_2}{N}\, Y \ge \frac{s_2}{N}\, Y_2 = \beta_2.
\]

The equality cases are verified by direct computation.
\end{proof}

\begin{proposition} \label{Proposition PS}
If $\seq{u_j}$ is a \PS{c} sequence of $E_\mu$, then a renamed subsequence converges weakly to a critical point $u$ of $E_\mu$, and either $u_j \to u$ in $W_\rad$ or
\[
c \ge E_\mu(u) + \beta^\ast \ge \beta^\ast.
\]
In particular, $E_\mu$ satisfies the {\em \PS{c}} condition for all $c < \beta^\ast$.
\end{proposition}

\begin{proof}
By Lemmas \ref{Lemma bounded BN} and \ref{Lemma weak limit}, a renamed subsequence of $\seq{u_j}$ converges weakly and a.e.\! to a critical point $u$ of $E_\mu$. Set $v_j = u_j - u$. By the Brezis-Lieb lemma \cite{MR699419}, applied to the functions $(u_j(x) - u_j(y))\, |x - y|^{- (N + s_i p_i)/p_i}$ in $L^{p_i}(\R^{2N})$ and to $u_j$ in $L^{p_i^\ast}(\R^N)$,
\[
X(u_j) = X(u) + X(v_j) + \o(1), \quad Y(u_j) = Y(u) + Y(v_j) + \o(1), \quad Z_i(u_j) = Z_i(u) + Z_i(v_j) + \o(1),
\]
and $\int_{\R^N} |u_j|^q\, dx \to \int_{\R^N} |u|^q\, dx$ by Theorem \ref{Theorem embedding} \ref{Theorem embedding.ii}. Since $E_\mu'(u_j)\, u_j \to 0$ by the boundedness of $\seq{u_j}$ and $E_\mu'(u)\, u = 0$, it follows that
\begin{gather}
\frac{1}{p_1}\, X(v_j) + \frac{1}{p_2}\, Y(v_j) - \frac{\mu_1}{p_1^\ast}\, Z_1(v_j) - \frac{\mu_2}{p_2^\ast}\, Z_2(v_j) = c - E_\mu(u) + \o(1), \label{E(v_j)}\\[7.5pt]
X(v_j) + Y(v_j) - \mu_1 Z_1(v_j) - \mu_2\, Z_2(v_j) = \o(1). \label{E'(v_j)}
\end{gather}
Suppose $u_j \not\to u$. Then, for a renamed subsequence, $\norm{v_j} \ge \delta$ for some $\delta > 0$ and
\[
X(v_j) \to X, \qquad Y(v_j) \to Y, \qquad \mu_1 Z_1(v_j) \to Z, \qquad \mu_2\, Z_2(v_j) \to W
\]
for some $X, Y, Z, W \ge 0$ with $(X,Y) \ne (0,0)$. Passing to the limit in \eqref{E(v_j)}, \eqref{E'(v_j)}, and \eqref{Sobolev Si} shows that $(X,Y,Z,W)$ satisfies \eqref{XYZW} and $I = c - E_\mu(u)$, so $c - E_\mu(u) \ge \beta^\ast$ by Lemma \ref{Lemma threshold}. Finally, $E_\mu(u) \ge 0$ by Lemma \ref{Lemma bounded BN}.
\end{proof}

\subsubsection{Proof of Theorem \ref{Theorem BN}}

We apply Theorem \ref{Theorem pseudo-index} to the even functional $E_\mu$ on $W_\rad$ with $c^\ast = \beta^\ast$; the \PS{c} condition for $c \in (0,\beta^\ast)$ holds by Proposition \ref{Proposition PS}.

\emph{The origin is a strict local minimizer.} For $\norm{u} \le 1$ we have $[u]_i \le 1$, so $[u]_2^{p_2} \ge [u]_2^{p_1}$ and
\[
\frac{1}{p_1}\, [u]_1^{p_1} + \frac{1}{p_2}\, [u]_2^{p_2} \ge \frac{1}{p_1} \left([u]_1^{p_1} + [u]_2^{p_1}\right) \ge \frac{1}{2^{p_1 - 1}\, p_1} \norm{u}^{p_1},
\]
while $\int_{\R^N} |u|^r\, dx \le C \norm{u}^r$ for $r \in \set{q,p_1^\ast,p_2^\ast}$ by Theorem \ref{Theorem embedding} \ref{Theorem embedding.i}. Since $q, p_1^\ast, p_2^\ast > p_1$, it follows that $E_\mu(u) \ge \norm{u}^{p_1} \left(2^{1 - p_1}/p_1 + \o(1)\right)$ as $u \to 0$.

\emph{The sets $C$, $A$, and $X$.} Let $C = \set{u \in V_m : \norm{u} = 1}$. Since $\norm{\cdot}$ is a norm on $V_m$, $C$ is the boundary of a bounded symmetric convex neighborhood of the origin in $V_m$, and hence is odd-homeomorphic to the unit sphere $S^{m-1}$; so $C$ is compact, symmetric, and $i(C) = m$. For $u \in C$ and $\tau \ge 0$, discarding the critical terms and using $[u]_i \le 1$ and \eqref{delta-m},
\begin{equation} \label{cone}
E_\mu(\tau u) \le \frac{\tau^{p_1}}{p_1} + \frac{\tau^{p_2}}{p_2} - \frac{\mu\, \delta_m}{q}\, \tau^q.
\end{equation}
Since $q > p_1 > p_2$, the right-hand side tends to $- \infty$ as $\tau \to \infty$; fix $R > 0$ so large that $\sup_{u \in A} E_\mu(u) \le 0$, where $A = \set{Ru : u \in C}$. Let $X = \set{tu : u \in A,\, 0 \le t \le 1} = \set{\tau u : u \in C,\, 0 \le \tau \le R}$. The condition $\mu > \nu_m$ in \eqref{nu-m} implies
\[
\frac{\tau^{p_1}}{p_1} + \frac{\tau^{p_2}}{p_2} - \frac{\mu\, \delta_m}{q}\, \tau^q < \beta^\ast \quad \forall \tau > 0,
\]
and the left-hand side is $0 < \beta^\ast$ at $\tau = 0$. Since $[0,R]$ is compact, \eqref{cone} gives $\sup_{u \in X} E_\mu(u) < \beta^\ast$. Theorem \ref{Theorem pseudo-index} now gives $m$ distinct pairs of critical points of $E_\mu$ at levels in $(0,\beta^\ast)$, which are nontrivial and, by Lemma \ref{Lemma symmetric criticality}, weak solutions of problem \eqref{BN}. Finally, $\nu_m < \infty$ for each $m$, so the number of solutions goes to infinity as $\mu \to \infty$. \qed

\begin{remark}
The proof of Theorem \ref{Theorem BN} uses only the following properties of the term $\mu\, |u|^{q-2}\, u$: it is compact on $W_\rad$, and its exponent exceeds $p_1$. Therefore, the same argument applies to problem \eqref{BN} with $\mu\, |u|^{q-2}\, u$ replaced by $\mu\, |x|^{-b}\, |u|^{q-2}\, u$ for any $(q,b) \in \Tc$ with $q > p_1$, with $\delta_m$ defined using $L$ in place of $\int_{\R^N} |u|^q\, dx$ and with $\vartheta = \min \set{q,p_1^\ast}$ in place of $p_1^\ast$ in Lemma \ref{Lemma bounded BN}.
\end{remark}

\subsection{Concluding remarks and open problems} \label{Section open problems}

\begin{enumroman}
\item {\em The region where \eqref{C} fails.} Theorem \ref{Theorem 2} $(ii)$ requires \eqref{C} only to verify $(f_4)$ via Theorem \ref{Theorem Pohozaev}, while Theorem \ref{Theorem 1} holds without \eqref{C}, which is needed there only to identify the eigenvalues $\widehat{\lambda}_k(q)$ with the minimax values $\lambda_k(q)$. Condition \eqref{C} reflects a symmetric H\"{o}lder splitting in Step 3 of the proof of Theorem \ref{Theorem Pohozaev}, and optimizing this splitting over $\H$ may enlarge the region where the identity is known.

Outside this region, the identity holds for solutions of class $C^{1,1}(\R^N)$ by an adaptation of the method of Ambrosio \cite{MR4761133}. It is an open problem whether every finite-energy weak radial solution of equation \eqref{15} with $(q,b) \in \H$ satisfies the Poho\v{z}aev identity \eqref{Pohozaev intro} without assuming \eqref{C}, in particular in the unweighted case $b = 0$.

A positive answer would remove condition \eqref{C} from Theorem \ref{Theorem 2} $(ii)$ and Theorem \ref{Theorem nonexistence}, and
would imply that $\widehat{\lambda}_k(q)=\lambda_k(q)$ for all $k$ in Theorem \ref{Theorem 1}.

\item {\em The subscaled region.} Theorem \ref{Theorem embedding} shows that the subscaled triangle $\Tsub$ lies in the range of the embeddings. Results for subscaled and concave-convex type problems with $(q,b) \in \Tsub$, in the spirit of \cite{PeWi}, will be considered elsewhere.
\item {\em The Sobolev points.} The existence or nonexistence of nontrivial radial solutions of the unweighted pure critical power equations remains open (see the remark following Theorem \ref{Theorem nonexistence}).
\end{enumroman}

\phantomsection
\addcontentsline{toc}{section}{Acknowledgments}
\section*{Acknowledgments}

\indent \indent KP was partially supported by Simons Foundation Travel Award 962241.

CS is a member of the Australian Mathematical Society (AustMS) and of INdAM-GNAMPA.

CS acknowledges the support of the Juan de la Cierva Fellowship (grant number JDC2023-050365-I).

\def\cprime{$''$}
\phantomsection


\addcontentsline{toc}{section}{References}
\begin{thebibliography}{10}

\bibitem{MR0370183}
Antonio Ambrosetti and Paul~H. Rabinowitz.
\newblock Dual variational methods in critical point theory and applications.
\newblock {\em J. Functional Analysis}, 14:349--381, 1973.

\bibitem{MR4761133}
Vincenzo Ambrosio.
\newblock On the {P}ohozaev identity for the fractional {$p$}-{L}aplacian
  operator in {$\R^N$}.
\newblock {\em Bull. Lond. Math. Soc.}, 56(6):1999--2013, 2024.

\bibitem{MR4597987}
G.~C. Anthal, J.~Giacomoni, and K.~Sreenadh.
\newblock A {C}hoquard type equation involving mixed local and nonlocal
  operators.
\newblock {\em J. Math. Anal. Appl.}, 527(2):Paper No. 127440, 27, 2023.

\bibitem{AnGa25}
Gurdev~Chand Anthal and Prashanta Garain.
\newblock Pohozaev-type identities for classes of quasilinear elliptic local
  and nonlocal equations and systems, with applications.
\newblock arXiv:2506.08667, 2025.

\bibitem{MR4785698}
Qiqi Bao.
\newblock Local-nonlocal {S}chr\"odinger equation with critical exponent: the
  zero mass case.
\newblock {\em Bull. Malays. Math. Sci. Soc.}, 47(5):Paper No. 163, 28, 2024.

\bibitem{MR3910033}
Mousomi Bhakta and Debangana Mukherjee.
\newblock Multiplicity results for {$(p,q)$} fractional elliptic equations
  involving critical nonlinearities.
\newblock {\em Adv. Differential Equations}, 24(3-4):185--228, 2019.

\bibitem{MR4992640}
Anup Biswas.
\newblock The {P}ohozaev identity for mixed local-nonlocal operators.
\newblock {\em J. Math. Anal. Appl.}, 557(1):Paper No. 130270, 19, 2026.

\bibitem{BiChDa25}
Nirjan Biswas, Souptik Chakraborty, and Paramananda Das.
\newblock Mixed local nonlocal equations with critical nonlinearity on
  $\mathbb{R}^n$: non-existence, existence, and multiplicity of positive
  solutions.
\newblock arXiv:2510.15507, 2025.

\bibitem{MR699419}
Ha{\"{\i}}m Br{\'e}zis and Elliott Lieb.
\newblock A relation between pointwise convergence of functions and convergence
  of functionals.
\newblock {\em Proc. Amer. Math. Soc.}, 88(3):486--490, 1983.

\bibitem{MR3007900}
X.~Chang and Z.-Q. Wang.
\newblock Ground state of scalar field equations involving a fractional
  {L}aplacian with general nonlinearity.
\newblock {\em Nonlinearity}, 26(2):479--494, 2013.

\bibitem{MR4626463}
Lassaad Chergui, Tianxiang Gou, and Hichem Hajaiej.
\newblock Existence and dynamics of normalized solutions to nonlinear
  {S}chr\"odinger equations with mixed fractional {L}aplacians.
\newblock {\em Calc. Var. Partial Differential Equations}, 62(7):Paper No. 208,
  45, 2023.

\bibitem{MR2018671}
Marco Degiovanni, Alessandro Musesti, and Marco Squassina.
\newblock On the regularity of solutions in the {P}ucci-{S}errin identity.
\newblock {\em Calc. Var. Partial Differential Equations}, 18(3):317--334,
  2003.

\bibitem{MR2944369}
Eleonora Di~Nezza, Giampiero Palatucci, and Enrico Valdinoci.
\newblock Hitchhiker's guide to the fractional {S}obolev spaces.
\newblock {\em Bull. Sci. Math.}, 136(5):521--573, 2012.

\bibitem{MR4736013}
Serena Dipierro, Kanishka Perera, Caterina Sportelli, and Enrico Valdinoci.
\newblock An existence theory for superposition operators of mixed order
  subject to jumping nonlinearities.
\newblock {\em Nonlinearity}, 37(5):Paper No. 055018, 27, 2024.

\bibitem{MR4942302}
Serena Dipierro, Kanishka Perera, Caterina Sportelli, and Enrico Valdinoci.
\newblock An existence theory for nonlinear superposition operators of mixed
  fractional order.
\newblock {\em Commun. Contemp. Math.}, 27(8):Paper No. 2550005, 39, 2025.

\bibitem{MR4821750}
Serena Dipierro, Edoardo Proietti~Lippi, Caterina Sportelli, and Enrico
  Valdinoci.
\newblock A general theory for the {$(s, p)$}-superposition of nonlinear
  fractional operators.
\newblock {\em Nonlinear Anal. Real World Appl.}, 82:Paper No. 104251, 24,
  2025.

\bibitem{DPSV26}
Serena Dipierro, Edoardo Proietti~Lippi, Caterina Sportelli, and Enrico
  Valdinoci.
\newblock Optimal embedding results for fractional {S}obolev spaces.
\newblock {\em Math. Z.}, 313:Paper No. 49, 40, 2026.

\bibitem{MR4685066}
Sidy~Moctar Djitte, Mouhamed~Moustapha Fall, and Tobias Weth.
\newblock A generalized fractional {P}ohozaev identity and applications.
\newblock {\em Adv. Calc. Var.}, 17(1):237--253, 2024.

\bibitem{MR0478189}
Edward~R. Fadell and Paul~H. Rabinowitz.
\newblock Generalized cohomological index theories for {L}ie group actions with
  an application to bifurcation questions for {H}amiltonian systems.
\newblock {\em Invent. Math.}, 45(2):139--174, 1978.

\bibitem{MR4142328}
Divya Goel, Deepak Kumar, and K.~Sreenadh.
\newblock Regularity and multiplicity results for fractional
  {$(p,q)$}-{L}aplacian equations.
\newblock {\em Commun. Contemp. Math.}, 22(8):1950065, 37, 2020.

\bibitem{MR1009077}
Mohammed Guedda and Laurent V{\'e}ron.
\newblock Quasilinear elliptic equations involving critical {S}obolev
  exponents.
\newblock {\em Nonlinear Anal.}, 13(8):879--902, 1989.

\bibitem{MR2784366}
Jiaquan Liu and Xiangqing Liu.
\newblock On the eigenvalue problem for the {$p$}-{L}aplacian operator in
  {$\R^N$}.
\newblock {\em J. Math. Anal. Appl.}, 379(2):861--869, 2011.

\bibitem{MR4445936}
Tingjian Luo and Hichem Hajaiej.
\newblock Normalized solutions for a class of scalar field equations involving
  mixed fractional {L}aplacians.
\newblock {\em Adv. Nonlinear Stud.}, 22(1):228--247, 2022.

\bibitem{MR5043800}
Carlo Mercuri and Kanishka Perera.
\newblock Variational methods for scaled functionals with applications to the
  {S}chr\"odinger-{P}oisson-{S}later equation.
\newblock {\em J. Math. Pures Appl. (9)}, 212:Paper No. 103885, 70, 2026.

\bibitem{MR547524}
Richard~S. Palais.
\newblock The principle of symmetric criticality.
\newblock {\em Comm. Math. Phys.}, 69(1):19--30, 1979.

\bibitem{Pe25}
Kanishka Perera.
\newblock A unifying zero-mass equation.
\newblock preprint, \href{https://arxiv.org/abs/2609.26585}{\tt
  arXiv:2609.26585}.

\bibitem{MR4999814}
Kanishka Perera.
\newblock Abstract multiplicity theorems and applications to critical growth
  problems.
\newblock {\em J. Anal. Math.}, 157(1):211--223, 2025.

\bibitem{KPCS}
Kanishka Perera and Caterina Sportelli.
\newblock A multiplicity result for critical elliptic problems involving
  differences of local and nonlocal operators.
\newblock {\em Topol. Methods Nonlinear Anal.}, 63(2):717--731, 2024.

\bibitem{PeWi}
Kanishka Perera and Patrick Winkert.
\newblock A scaling approach to double phase equations in $\mathbb{R}^N$.
\newblock preprint.

\bibitem{MR192184}
S.~I. Poho{\v{z}}aev.
\newblock On the eigenfunctions of the equation {$\Delta u+\lambda f(u)=0$}.
\newblock {\em Dokl. Akad. Nauk SSSR}, 165:36--39, 1965.

\bibitem{MR855181}
Patrizia Pucci and James Serrin.
\newblock A general variational identity.
\newblock {\em Indiana Univ. Math. J.}, 35(3):681--703, 1986.

\bibitem{MR845785}
Paul~H. Rabinowitz.
\newblock {\em Minimax methods in critical point theory with applications to
  differential equations}, volume~65 of {\em CBMS Regional Conference Series in
  Mathematics}.
\newblock Conference Board of the Mathematical Sciences, Washington, DC; by the
  American Mathematical Society, Providence, RI, 1986.

\bibitem{MR3211861}
Xavier Ros-Oton and Joaquim Serra.
\newblock The {P}ohozaev identity for the fractional {L}aplacian.
\newblock {\em Arch. Ration. Mech. Anal.}, 213(2):587--628, 2014.

\end{thebibliography}
\end{document}